\documentclass[12pt,leqno]{amsart}
\usepackage{amsmath,amsthm,amsfonts,amsbsy,multicol,fancybox}
\usepackage[utf8]{inputenc}
\usepackage[T1]{fontenc}

\usepackage{caption}\usepackage{subcaption}
\usepackage{graphicx}\usepackage{amssymb,latexsym}
\usepackage{epstopdf}
\usepackage{color}
\DeclareGraphicsExtensions{.pdf, .jpg, .pnp, .bmp, .fig}
\usepackage{tikz}
\usepackage[all]{xy}
\input xy
\xyoption{all}

\usepackage[hidelinks]{hyperref}

\usepackage{thmtools}

\newcounter{argument}
\newenvironment{argument}[1][\medskip]{%
\refstepcounter{argument}
\par\medskip
\noindent\phantomsection
\textbf{#1~\thesection.\arabic{argument}\,\,}\rmfamily\em}{\hspace{\fill}$\Box$\par\smallskip\noindent}

\newcommand{\bass}{\begin{argument}[Assumption]}\newcommand{\eass}{\end{argument}}
\newcommand{\bth}{\begin{argument}[Theorem]} \newcommand{\ethe}{\end{argument}}
\newcommand{\bre}{\begin{argument}[Remark]}      \newcommand{\ere}{\end{argument}}
\newcommand{\ble}{\begin{argument}[Lemma]}       \newcommand{\ele}{\end{argument}}
\newcommand{\bde}{\begin{argument}[Definition]}   \newcommand{\ede}{\end{argument}}
\newcommand{\bco}{\begin{argument}[Corollary]}     \newcommand{\eco}{\end{argument}}
\newcommand{\bpr}{\begin{argument}[Proposition]}  \newcommand{\epr}{\end{argument}}
\newcommand{\bexam}{\begin{argument}[Example]}\newcommand{\eexam}{\end{argument}}

\newcommand{\bpf}{\begin{proof}}\newcommand{\epf}{\end{proof}}
\newcommand{\barr}{\begin{array}}\newcommand{\earr}{\end{array}}
\newcommand{\beao}{\begin{eqnarray*}}\newcommand{\eeao}{\end{eqnarray*}\noindent}
\newcommand{\beam}{\begin{eqnarray}}\newcommand{\eeam}{\end{eqnarray}\noindent}
\newcommand{\beqq}{\begin{equation}}\newcommand{\eeqq}{\end{equation}\noindent}

 \newcommand{\un}{\underbrace}
\newcommand{\wt}{\widetilde}

\newcommand{\nto}{n\to\infty}  
\newcommand{\kto}{k\to\infty} \newcommand{\mto}{m\to\infty}

\newcommand{\tto}{t\to\infty}

\newcommand{\del}{\delta}
\newcommand{\D}{\Delta}
  \newcommand{\ep}{\epsilon}

\newcommand{\lam}{\lambda} 

\newcommand{\w}{\omega} \newcommand{\W}{\Omega}

\newcommand{\bba}{{\mathcal A}}
\newcommand{\bbb}{{\mathcal B}}
\newcommand{\bbc}{{\mathcal C}} 

\newcommand{\bfE}{{\mathbb E}}\newcommand{\bbE}{{\mathcal E}} 
\newcommand{\bbf}{{\mathcal F}}

\newcommand{\bbi}{{\mathbb I}}

\newcommand{\bbl}{{\mathcal L}}

 \newcommand{\bbN}{{\mathbb N}}

\newcommand{\bfP}{{\mathbb P}}
                              
 \newcommand{\bbR}{{\mathbb R}}

\newcommand{\bbt}{{\mathcal T}}

\graphicspath{{./../Images/}}

\begin{document}

\title[CIR/CEV type delay models with jump]{An explicit positivity preserving numerical scheme for CIR/CEV type delay models with jump}

\author[I. S. Stamatiou]{I. S. Stamatiou}
\email{joniou@gmail.com, ioannis.stamatiou@ouc.ac.cy, istamatiou@uniwa.gr}

\begin{abstract}
We consider mean-reverting CIR/CEV processes with delay and jumps used as models on the financial markets. These processes are solutions of stochastic differential equations with jumps, which have no explicit solutions. We prove the non-negativity property of the solution of the above models and propose an explicit positivity preserving numerical scheme,using the semi-discrete method, that converges in the strong sense to the exact solution. We also make some minimal numerical experiments to illustrate the proposed method.
\end{abstract}

\date\today

\maketitle
\tableofcontents\listoffigures\listoftables

\section{Introduction}\label{CIRCEVdj:sec:intro}
\setcounter{equation}{0}

Throughout, let $T>0$ and $(\Omega, \bbf, \{\bbf_t\}_{0\leq t\leq T}, \bfP)$ be a complete probability space, meaning that the filtration $ \{\bbf_t\}_{0\leq t\leq T} $ satisfies the usual conditions, i.e. is right continuous and $\bbf_0$ includes all $\bfP-$null sets. Let $W_{t,\w}:[0,T]\times\W\rightarrow\bbR$ be a one-dimensional Wiener process adapted to the filtration $\{\bbf_t\}_{0\leq t\leq  T}$ and $\wt{N}(t)=N(t)-\lam t$ a compensated Poisson process with intensity $\lam>0$ independent of  $W_t.$ Consider the following stochastic delay differential equation (SDDE) with jump,
\beqq\label{CIRCEVdj:eq:sddej}
x_t=\left\{ \barr{ll}
\xi_0 + \int_{0}^{t}(k_1 - k_2x_{s^-})ds + \int_{0}^{t}k_3b(x_{s-\tau})x_{s^-}^\alpha dW_s + \int_{0}^{t}g(x_{s^-}) d\wt{N}_s, & t\in [0,T],\\
&\\
\xi(t),& t\in[-\tau,0],
\earr \right.
\eeqq
where $x_{s^-}=\lim_{r\uparrow s}x_r,$ the coefficient $b\in \bbc(\bbR_+,\bbR_+),$\footnote{$\bbc(A,B)$ the space of continuous functions $\phi:A\mapsto B$ with norm $\|\phi\|=\sup_{u\in A}\phi(u)$} $g:\bbR\mapsto \bbR$ is the jump coefficient assumed deterministic for simplicity, the function $\xi\in \bbc([-\tau,0],(0,\infty))$ and $\tau>0$ is a positive constant which represents the delay. The quantities $k_i$ are positive constants and the real number $\alpha\in[1/2,1).$ The diffusion coefficient does not satisfy standard assumptions, i.e. linear growth conditions, therefore classical results on the existence and uniqueness of solution and boundness of the moments as in \cite{mao:1994} and \cite{kuchler_platen:2000} are not applicable here. Model (\ref{CIRCEVdj:eq:sddej}) includes many models used in mathematical finance describing quantities such as interest rates and volatilities, see \cite{platen:2010} and Table \ref{CIRCEVdj:table:t1}.\footnotesize

\begin{table}[h!]
	\begin{center}
		\caption{Models included in general model (\ref{CIRCEVdj:eq:sddej}).}
		\label{CIRCEVdj:table:t1}
		\begin{tabular}{rrr|l} 
			Jump Coefficient & Delay coefficient & Diffusion exponent & Model/Ref\\
			$g(x)$ & $b(x)$ & $\alpha$ &\\
			\hline
			0 & 1 & 1/2 & CIR, \cite{cox_et_al:1985}\\
			0 & $x^\gamma$ & 1/2&delay CIR, \cite{wu_etal:2009},\cite{fatemion_etal:2017}\\
			$\del x$ & 1& 1/2&CIR with jump, \cite{wu_etal:2008}, \cite{wu:2014}\\
			$\del x$ & $x^\gamma$ & 1/2& delay CIR with jump, \cite{jiang_etal:2011}, \cite{fatemion_etal:2018}\\
			\hline
			0 & 1 & (1/2,1) & CEV, \cite{cox:1975}\\
			0 & $x^\gamma$ & (1/2,1)&delay CEV,\\
			$\del x$ & 1& (1/2,1)& CEV with jump, \cite{beliaeva:2012}\\
			$\del x, \del\sin x, \del x/(1+x)$ & 1& [1/2,1)& CIR/CEV with jump, \cite{yang_etal:2018}
		\end{tabular}
	\end{center}
\end{table}\normalsize

Numerical approximations of SDDEs with jumps like (\ref{CIRCEVdj:eq:sddej}) are necessary for simulations of the paths $x_t(\w),$ or for approximation of functions of $x(T)$ or functionals of the form $\bfE F(x),$ where $F:\bbc([0,T],\bbR)\mapsto\bbR,$ so as to obtain the expected payoff of an option. The numerical analysis of jump models become more computationally complex proportional to the jump intensity \cite{bruti_liberati_platen:2007}.\\
We assume the following conditions for the delay coefficient $b$ and the jump coefficient $g.$\\
\textbf{Assumption A}  
The delay diffusion coefficient $b(\cdot)$ is $\gamma$-H\"{o}lder continuous where $\gamma>0$ i.e.
$$sup_{z<x}\frac{|b(x) - b(z)|}{(x-z)^\gamma}:=C_\gamma<\infty,$$
and the jump coefficient $g(\cdot)$ satisfies for a constant $L>0$ the following relation $\forall x\in\bbR$  
$$\Big(|g^{\prime}(x)|<L\leq1,  \hbox{ and } g(0)=0 \hbox{ OR } g(x)>0 \hbox{ and } |g^{\prime}(x)|<L.\Big)$$
Moreover assume $\|\xi\|<A_\xi$ for a constant $A_\xi$.\\
Let  $\bbt$ be a jump-adapted time partition of the interval $[-\tau,T],$ see Section \ref{CIRCEVdj:sec:main} for details. We introduce the following numerical method for the approximation of the solution of (\ref{CIRCEVdj:eq:sddej}), for the jump-extended CIR/CEV model with delay, which we call Jump Adapted Semi-Discrete method (JASDM) with $y_t=\xi(t)$  when $t\in[-\tau,0]$ and for $k=0, 1, \ldots, n_T-1,$
\beqq\label{CIRCEVdj:eq:JASD}
\left\{\barr{ll}
y_{t_{{k+1}^-}} = (z_{t_{k+1}})^2,\\
y_{t_{k+1}} = y_{t_{{k+1}^-}} + g(y_{t_{{k+1}^-}})\D \wt{N}_k,\\
\earr \right.
\eeqq
where \footnotesize
\beao 
z_{t} &=&  \sqrt{y_{t_k}\left(1- \frac{k_2\D_k}{1+k_2\theta\D_k}\right) + \frac{k_1\D_k}{1+k_2\theta\D_k} - \frac{(k_3)^2}{4(1+k_2\theta\D_k)^2 }\frac{b^2(y_{t_{k-\tau}})}{(1+b(y_{t_{k-\tau}})\D_k^{m})^2}(y_{t_k})^{2\alpha-1}\D_k}\\
&& +\frac{k_3}{2(1+k_2\theta\D_k)} \frac{b(y_{t_{k-\tau}})}{1+b(y_{t_{k-\tau}})\D_k^{m}}(y_{t_k})^{\alpha -\frac{1}{2}}(W_t-W_{t_k})
\eeao
\normalsize
and
$\D_k=t_{k+1}-t_{k}, \D \wt{N}_k:=\wt{N}(t_{k+1})-\wt{N}(t_{k})=\D N_k-\lam\D_k$ and $\theta\in[0,1]$ represents the level of implicitness. The positive constant $m$ is chosen to be equal to $1/4,$ see proof of Proposition \ref{CIRCEVdj-propo:L2auxiliaryConv}. The increment $\D N_k=1$ when $t_{k+1}$ is a jump time and zero otherwise. For the derivation of (\ref{CIRCEVdj:eq:JASD}) briefly saying we apply an analogue of the semi-discrete method \cite{halidias_stamatiou:2015} (and \cite{halidias:2015} for the case $\alpha=1/2$) for the numerical approximation of the SDDE (\ref{CIRCEVdj:eq:sddej}) without jump, that is in between jump times with the appropriate modification when the jump times occur, see Section \ref{CIRCEVdj:sec:main} for details.   

The numerical scheme (\ref{CIRCEVdj:eq:JASD}) is positivity preserving by construction and converges to the true solution $x_t$ of model (\ref{CIRCEVdj:eq:sddej}) and this is stated in our main result, Theorem  \ref{CIRCEVdj:thm:strong_convSDDEjump}.

A basic feature of the semi-discrete method, originally proposed in \cite{halidias:2012}, is that it is explicit,  converges strongly in the mean square sense to the exact solution of the original SDE, preserves positivity (see application in Wright-Fisher model \cite{stamatiou:2018}),  and it does not explode in some super-linear problems \cite{halidias:2014},  \cite{halidias_stamatiou:2016}. The purpose of this paper is on the one hand to  propose the general model (\ref{CIRCEVdj:eq:sddej}) and to generalize the semi-discrete method to include SDDEs with jumps that arise in the area of financial mathematics in valuating options. 

There has been a lot of effort in numerical approximations of models of type (\ref{CIRCEVdj:eq:sddej}) mainly in cases where no jumps are included, and we simply refer to Table \ref{CIRCEVdj:table:t1} and references therein. The main difficulty in producing strongly converging numerical schemes, as already stated, is the non-Lipschitz nature of the coefficients. The recent focus on positivity preserving methods for these problems can be stated in   \cite{yang_etal:2018}, \cite{fatemion_etal:2017}.

In Section \ref{CIRCEVdj:sec:main} we collect all the main results. First, we investigate the nonnegative solutions of (\ref{CIRCEVdj:eq:sddej}) and the mean-reversion property. Then we examine properties of the proposed numerical scheme (\ref{CIRCEVdj:eq:JASD}), such as strong convergence, boundness of moments and the mean-reversion property. Section \ref{CIRCEVdj:sec:numerics} gives some numerical examples confirming theoretical results. 
Section \ref{CIRCEVdj:sec:proofs} is devoted to the proofs of all the results and in the last subsection we provide some concluding remarks.

\section{Setting and main results}\label{CIRCEVdj:sec:main}

Our first goal, is to provide the well-posedness of (\ref{CIRCEVdj:eq:sddej}), that is the existence and uniqueness of solution of the CIR/CEV model with delay, as well as the positivity of the solution process.  Instead of proving existence, uniqueness and non-negativity of solution of (\ref{CIRCEVdj:eq:sddej}) directly, we treat the following SDDE with jump 
\footnotesize
\beqq\label{CIRCEVdj:eq:sddejABS}
x_t=\left\{ \barr{ll}
\xi_0 + \int_{0}^{t}(k_1 - k_2x_{s^-})ds + \int_{0}^{t}k_3b(|x_{s-\tau}|)|x_{s^-}|^\alpha dW_s + \int_{0}^{t}g(x_{s^-}) d\wt{N}_s, & t\in [0,T],\\
&\\
\xi(t),& t\in[-\tau,0].
\earr \right.
\eeqq
\normalsize
We first show that (\ref{CIRCEVdj:eq:sddejABS}) has finite $p$-moments.

\ble[Moment bounds for (\ref{CIRCEVdj:eq:sddejABS})]
\label{CIRCEVdj:lem:moment_bound}
Let $x_t$ be the solution of (\ref{CIRCEVdj:eq:sddejABS}). Then for any $p>0$ there exists a positive constant $A$ such that
\beqq \label{CIRCEVdj:eq:moment_bound}
\bfE\sup_{-\tau\leq t\leq T}|x_t|^p\leq A.
\eeqq
\ele

We then examine the uniqueness and positivity of (\ref{CIRCEVdj:eq:sddejABS}).

\bpr[Uniqueness and positivity of (\ref{CIRCEVdj:eq:sddejABS})]
\label{CIRCEVdj:prop:uniq_posABS}
There exists a unique and nonnegative solution for equation (\ref{CIRCEVdj:eq:sddejABS}), i.e. $x_t>0$ almost surely (a.s.).
\epr
Using the previous results we easily reach the following theorem which is our first main result. We also provide the proof since it is short.

\bth[Existence, Uniqueness and positivity of (\ref{CIRCEVdj:eq:sddej})] \label{CIRCEVdj:thm:uniq_pos}
There exist a  unique positive solution for the jump-extended CIR/CEV model (\ref{CIRCEVdj:eq:sddej}) and the $p$-th moment of the solution is bounded. 
\ethe

\bpf[Proof of Theorem \ref{CIRCEVdj:thm:uniq_pos}]
It is an immediate consequence of Lemma \ref{CIRCEVdj:lem:moment_bound} and Proposition \ref{CIRCEVdj:prop:uniq_posABS} since when $x$ positive $|x|=x.$
\epf

We also give inverse moment bounds of $(x_t)$ which will be used in the proof of our main convergence results, Theorems \ref{CIRCEVdj:thm:strong_convSDDE} and \ref{CIRCEVdj:thm:strong_convSDDEjump}.

\ble[Inverse Moment bounds for (\ref{CIRCEVdj:eq:sddej})]
\label{CIRCEVdj:lem:inv_moment_bound}
Let $x_t$ be the solution of (\ref{CIRCEVdj:eq:sddej}), then for any $p<0$ there exists a positive constant $\hat{A}$ such that
\beqq \label{CIRCEVdj:eq:inv_moment_bound}
\sup_{-\tau\leq t\leq T}\bfE(x_t)^p \leq \hat{A}.
\eeqq
\ele

We now examine the mean reversion property which holds for model (\ref{CIRCEVdj:eq:sddej}) with no delay and jump, see \cite[Theorem 2.2]{wu_etal:2008}.

\bth[Mean-reversion] \label{CIRCEVdj:thm:mean_reversion}
The jump-extended CIR/CEV model (\ref{CIRCEVdj:eq:sddej}) preserves the mean reversion property, i.e.
\beqq\label{CIRCEVdj:eq:mean_reversion}
\lim_{\tto}\bfE x_t = \frac{k_1}{k_2}.
\eeqq
\ethe

\bpf[Proof of Theorem \ref{CIRCEVdj:thm:mean_reversion}]
We take expectations on both sides of (\ref{CIRCEVdj:eq:sddej}) to get
$$
\bfE x_t =\bfE\xi_0 - \int_{0}^{t}(k_2\bfE x_{s^-} - k_1)ds,
$$ 
where we also used the boundeness of the moments of $(x_t)$ and the fact that $(\wt{N}_s)$ is a martingale. Now, let $\psi_t:=\bfE x_t.$ Then $\psi_t$ satisfies the following integral equation 
$$
\psi_t = \psi_0 - \int_{0}^{t}(k_2\psi_{s^-}- k_1)ds,
$$
or $(\psi_t)^{\prime}+ k_2\psi_t =k_1$ and multiplying by $e^{k_2t}$ each side of the equality implies $\psi_t=\frac{k_1}{k_2} + \left(\bfE \xi_0 -\frac{k_1}{k_2}\right)e^{-k_2t}$. We therefore verify (\ref{CIRCEVdj:eq:mean_reversion}). 
\epf

We rewrite the stochastic integral equation (\ref{CIRCEVdj:eq:sddej}), denoting by a superscript $g$ the dependence on the jump coefficient as 
$$
x_t^g=\left\{ \barr{ll}
\xi_0 + \int_{0}^{t}(k_1 - k_2x_{s^-})ds + \int_{0}^{t}k_3b(x_{s-\tau})x_{s^-}^\alpha dW_s + \int_{0}^{t}g(x_{s^-}) d\wt{N}_s, & t\in [0,T],\\
&\\
\xi(t),& t\in[-\tau,0].
\earr \right.
$$
and in a similar way we denote the semi-discrete scheme $y_t^g.$ Let the horizon $T$ be a multiple of $\tau,$ i.e. $T=N_0\tau,$ where $N_0\in\bbN$ and discretize the interval $[-\tau,T]$ with equidistant steps of size $\D=\tau/l$ for some $l\in\{2,3,\ldots\}$ and $t_n=n\D=n\tau/l,$ where $n=-l,-l+1,\ldots, N.$ Now, construct a jump-adapted time partition 
$$\bbt=\{-\tau=t_{-l}< t_{-l+1}<\ldots<0=t_0<t_1<\ldots<t_{N_T}=T\},$$
produced by a superposition of the jump times to the deterministic equidistant grid. Here,  $N_t=\max\{n: t_n\leq t\}$  for $t\leq T$ and $\bfP(t_{n+1}-t_n\leq\D)=1.$ The discretization $\bbt$ is path-dependent with maximum step-size $\D.$ For $t\in[t_k,t_{k+1})$ and $k\geq0$ it holds 
$$
\left\{ \barr{ll}
x_{t_{{k+1}^-}}=x_{t^-} + \int_{t}^{t_{k+1}}(k_1 - k_2x_{s^-})ds + \int_{t}^{t_{k+1}}k_3b(x_{s-\tau})x_{s^-}^\alpha dW_s,\\
&\\
x_{t_{k+1}} = x_{t_{{k+1}^-}} + g(x_{t_{{k+1}^-}}) \D\wt{N}_k,
\earr \right.
$$
where  $\D \wt{N}_k:= N((k+1)\D) - N(k\D) - \lam\D$ are the compensated Poisson increments. Between jump times (\ref{CIRCEVdj:eq:sddej}) evolves like a delay CIR/CEV process without jumps. Theorem \ref{CIRCEVdj:thm:uniq_pos} implies that also the delay CIR/CEV model with no jumps has a unique nonnegative solution, i.e. $x_t^0\geq0$ a.s.
Therefore, the idea is to approximate the process between jump  times with a positivity preserving numerical method and then take into account the jumps. In the following we will use the ideas in \cite{yang_etal:2018} where the mean-reverting CIR/CEV jump process was considered with no delay.\\
Let
\beam\nonumber
f_{\theta}(x,y)&=& \underbrace{k_1 - k_2(1-\theta)x - \frac{(k_3)^2}{4(1+k_2\theta\D)}\frac{b^2(z)}{(1+b(z)\D^{m})^2}x^{2\alpha-1}-k_2\theta y}_{f_1(x,y,z)}\\
\label{CIRCEVdj-eq:drift discetisation}&& + \underbrace{\frac{(k_3)^2}{4(1+k_2\theta\D)}\frac{b^2(z)}{(1+b(z)\D^{m})^2}x^{2\alpha-1}}_{f_2(x,z)}
\eeam
and
\beqq\label{CIRCEVdj-eq:diffusion discetisation}
g(x,y,z)=k_3 \frac{b(z)}{1+b(z)\D^{m}}x^{\alpha-\frac{1}{2}}\sqrt{y},
\eeqq
where $f_\theta(x,x,x)=k_1-k_2 x$ and $g(x,x,z)=k_3\frac{b(z)}{1+b(z)\D^{m}}x^a.$

Consider the following process
$$
y_t=y_{t_n} + f_1(y_{t_n}, y_t, y_{t_{n}-\tau})\cdot\D + \int_{t_n}^tf_2(y_{t_n}, y_{t_{n}-\tau})ds + \int_{t_n}^{t} \textup{sgn}(z_s)g(y_{t_n},y_s, y_{t_{n}-\tau})dW_s,
$$
with $y_0=\xi_0$ a.s. or more explicitly\footnotesize
\beam
\nonumber
&&y_t=y_{t_n} + \left(k_1 - k_2(1-\theta)y_{t_n} - \frac{(k_3)^2}{4(1+k_2\theta\D)}\frac{b^2(y_{t_{n}-\tau})}{(1+b(y_{t_{n}-\tau})\D^{m})^2}(y_{t_n})^{2\alpha-1}-k_2\theta y_t\right)\cdot\D\\
\nonumber&&+ \int_{t_n}^{t}\frac{(k_3)^2}{4(1+k_2\theta\D)}\frac{b^2(y_{t_{n}-\tau})}{(1+b(y_{t_{n}-\tau})\D^{m})^2}(y_{t_n})^{2\alpha-1} ds\\
&&\label{CIRCEVdj-eq:SDsde}
+ k_3\frac{b(y_{t_{n}-\tau})}{1+b(y_{t_{n}-\tau})\D^{m}}(y_{t_n})^{\alpha-\frac{1}{2}}\int_{t_n}^{t} \textup{sgn}(z_s)\sqrt{y_s}dW_s,
\eeam\normalsize
for $ t\in(t_n,t_{n+1}],$ where $\theta\in[0,1]$ represents the level of implicitness and
\beqq\label{CIRCEVdj-eq:SDsgn term}
z_t= \sqrt{y_n} +\frac{k_3}{2(1+k_2\theta\D)}\frac{b(y_{t_{n}-\tau})}{1+b(y_{t_{n}-\tau})\D^{m}}(y_{t_n})^{\alpha-\frac{1}{2}}(W_t-W_{t_n}),
\eeqq
with
\beqq\label{CIRCEVdj-eq:SDinside sgn term}
y_n:=y_{t_n}\left(1- \frac{k_2\D}{1+k_2\theta\D}\right) + \frac{k_1\D}{1+k_2\theta\D} - \frac{(k_3)^2}{4(1+k_2\theta\D)^2}\frac{b^2(y_{t_{n}-\tau})}{(1+b(y_{t_{n}-\tau})\D^{m})^2}(y_{t_n})^{2\alpha-1}\D.
\eeqq
Note that $y_n=J_n/(1+k_2\theta\D)$ where
\beao
J_n&=&\left(y_{t_n}(1-k_2(1-\theta)\D) + k_1\D - \frac{(k_3)^2}{4(1+k_2\theta\D)}\frac{b^2(y_{t_{n}-\tau})}{(1+b(y_{t_{n}-\tau})\D^{m})^2}(y_{t_n})^{2\alpha-1}\D\right)\bbi_{y_{t_n}\geq1}\\
&&+\left(y_{t_n}(1-k_2(1-\theta)\D) + k_1\D - \frac{(k_3)^2}{4(1+k_2\theta\D)}\frac{b^2(y_{t_{n}-\tau})}{(1+b(y_{t_{n}-\tau})\D^{m})^2}(y_{t_n})^{2\alpha-1}\D\right)\bbi_{y_{t_n}<1}\\
&\geq&y_{t_n}^{2\alpha-1}\left(1-k_2(1-\theta)\D - \frac{(k_3)^2}{4(1+k_2\theta\D)}\sqrt{\D}\right)\bbi_{y_{t_n}\geq1}\\
&&+\left(y_{t_n}(1-k_2(1-\theta)\D) - \frac{(k_3)^2}{4(1+k_2\theta\D)}\sqrt{\D}(y_{t_n})^{2\alpha-1}\right)\bbi_{\D^{1/4}<y_{t_n}<1}\\
&&+\left(y_{t_n}(1-k_2(1-\theta)\D) +k_1\D- \frac{(k_3)^2}{4(1+k_2\theta\D)}\sqrt{\D}(y_{t_n})^{2\alpha-1}\right)\bbi_{y_{t_n}\leq\D^{1/4}}\\
&\geq&1-k_2(1-\theta)\D - \frac{(k_3)^2}{4(1+k_2\theta\D)}\sqrt{\D}\\
&&+y_{t_n}\left(1-k_2(1-\theta)\D - \frac{(k_3)^2}{4(1+k_2\theta\D)}\D^{1/4}\right)\bbi_{\D^{1/4}<y_{t_n}<1}\\
&&+\left(y_{t_n}(1-k_2(1-\theta)\D) +k_1\D- \frac{(k_3)^2}{4(1+k_2\theta\D)}\D^{(2\alpha+1)/4}\right)\bbi_{y_{t_n}\leq\D^{1/4}},
\eeao
and $m\leq1/4$ to be specified later. In the above inequality, the first term is nonnegative when $\D<(k_2(1-\theta) + (k_3)^2/4)^{-2}\wedge(k_2(1-\theta))^{-1}$ the second term is nonnegative when $\D<(k_2(1-\theta) + (k_3)^2/4)^{-4}$ and for the last term we have
\beao
&&\left(y_{t_n}(1-k_2(1-\theta)\D) +k_1\D- \frac{(k_3)^2}{4(1+k_2\theta\D)}\D^{(2\alpha+1)/4}\right)\bbi_{y_{t_n}\leq\D^{1/4}}\\
&\geq&\left(y_{t_n}\frac{(k_3)^2}{4}  + k_1\D- \frac{(k_3)^2}{4(1+k_2\theta\D)}\sqrt{\D}\right)\bbi_{y_{t_n}\leq\D^{1/4}}\geq0,
\eeao
when $\D<\frac{1-(k_3)^2/4}{k_2(1-\theta)}.$

Furthermore, (\ref{CIRCEVdj-eq:SDsde}) has jumps at nodes $t_n.$ Solving for $y_t,$ we end up with the following explicit scheme
\beao
y_t&=&y_n + \int_{t_n}^{t}\frac{(k_3)^2}{4(1+k_2\theta\D)^2}\frac{b^2(y_{t_{n}-\tau})}{(1+b(y_{t_{n}-\tau})\D^{m})^2}(y_{t_n})^{2\alpha-1}ds\\
&&+ \frac{k_3}{1+k_2\theta\D}\frac{b(y_{t_{n}-\tau})}{1+b(y_{t_{n}-\tau})\D^{m}}(y_{t_n})^{\alpha-\frac{1}{2}}\int_{t_n}^{t} \textup{sgn}(z_s)\sqrt{y_s} dW_s,
\eeao 
with solution in each step given by \cite[(4.39), p.123]{kloeden_platen:1995}
$$
y_t=(z_t)^{2}.
$$

We remove the term $\textup{sgn}(z_s)$ from (\ref{CIRCEVdj-eq:SDsde}) by considering the process
$$ 
\wt{W}_t=\int_{0}^{t} \textup{sgn}(z_s)dW_s,
$$ 
which is a martingale with quadratic variation $< \wt{W}_t, \wt{W}_t>=t$ and thus a standard Brownian motion w.r.t. its own filtration, justified by L\'evy's theorem \cite[Th. 3.16, p.157]{karatzas_shreve:1988}. Therefore, the compact form of (\ref{CIRCEVdj-eq:SDsde}) becomes
\beam
\nonumber&&y_t=\xi_0 + \int_0^t\left(k_1 - k_2(1-\theta)y_{\hat{s}}-k_2\theta y_{\wt{s}}\right)ds\\
&&\nonumber + \int_t^{t_{n+1}}\left(k_1 - k_2(1-\theta)y_{t_n} - \frac{(k_3)^2}{4(1+k_2\theta\D)}\frac{b^2(y_{t_{n}-\tau})}{(1+b(y_{t_{n}-\tau})\D^{m})^2}(y_{t_n})^{2\alpha-1}-k_2\theta y_{t}\right)ds\\
&&\label{CIRCEVdj-eq:SDsdeComp}
+ k_3\int_{0}^{t} \frac{b(y_{\hat{s}-\tau})}{1+b(y_{\hat{s}-\tau})\D^{m}}(y_{\hat{s}})^{\alpha-\frac{1}{2}}\sqrt{y_s}d\wt{W}_s,
\eeam
for $t\in(t_n,t_{n+1}]$  where
$$
\hat{s}=t_j, s\in(t_j,t_{j+1}], \, j=0,\ldots,n, \quad
\wt{s}=\left\{ \barr{ll}  t_{j+1},  & \mbox{for } \, s\in[t_j,t_{j+1}], \\  t,  & \mbox{for }\,  s\in[t_n,t] \earr j=0,\ldots,n-1.\right.
$$
Consider also the process
\beqq\label{CIRCEVdj-eq:CEVunderNewWiener}
\wt{x}_t=\xi_0 + \int_0^t (k_1-k_2 \wt{x}_s)ds + \int_{0}^{t}k_3b(\wt{x}_{s-\tau})(\wt{x}_s)^\alpha d\wt{W}_s, \quad t\in [0,T],
\eeqq
with $\wt{x}_t=\xi(t), t\in[-\tau,0].$ The process $(x_t^0)$ and the process $(\wt{x}_t)$ of (\ref{CIRCEVdj-eq:CEVunderNewWiener}) have the same distribution. We show in the following that $\bfE\sup_{0\leq t\leq T}|y_t-\wt{x}_t|^2\rightarrow0$ as $\D\downarrow0$ thus $\bfE\sup_{0\leq t\leq T}|y_t-x_t|^2\rightarrow0$ as $\D\downarrow0.$ To simplify notation we write $\wt{W}, (\wt{x}_t)$ as $W, (x_t).$
We end up with the following explicit scheme
\beam\nonumber
y_t^0&=&y_n + \int_{t_n}^{t}\frac{(k_3)^2}{4(1+k_2\theta\D)^2}\frac{b^2(y_{t_{n}-\tau})}{(1+b(y_{t_{n}-\tau})\D^{m})^2}(y_{t_n})^{2\alpha-1}ds\\
&&\label{CIRCEVdj-eq:SDexplNW}
+ \frac{k_3}{1+k_2\theta\D}\frac{b(y_{t_{n}-\tau})}{1+b(y_{t_{n}-\tau})\D^{m}}(y_{t_n})^{\alpha-\frac{1}{2}}\int_{t_n}^{t} \sqrt{y_s} dW_s,
\eeam
where $y_n$ is as in (\ref{CIRCEVdj-eq:SDinside sgn term}).

The following result shows that the (JASDM) method preserves the mean-reversion property for appropriate step size $\D$ and implicitness parameter $\theta.$ 

\bth[Mean-reversion Property] \label{CIRCEVdj:thm:mean_reversionJASD}
The (JASDM) method (\ref{CIRCEVdj:eq:JASD}) preserves the mean reversion property in the following sense, 
\beqq\label{CIRCEVdj:eq:mean_reversionJASD}
\lim_{\nto}\bfE y_{t_n} \leq \frac{k_1}{k_2\theta}.
\eeqq
when $\D(1-\theta)<1/k_2$ and $0<\theta\leq1.$
\ethe
\textbf{Assumption B}  
Let $\D>0$ be such that \footnotesize $\D<\left(\frac{1}{k_2(1-\theta) + (k_3)^2/4}\right)^{2}\wedge \left(\frac{1}{k_2(1-\theta) + (k_3)^2/4)}\right)^{4}\wedge\frac{4-(k_3)^2}{4k_2(1-\theta)}$  \normalsize for $\theta\in[0,1].$ 
\\ 
In the following we show the strong convergence of the semi-discrete method to the true solution of model (\ref{CIRCEVdj:eq:sddej}) without jumps. 

\bth[Strong convergence when no jumps are included] \label{CIRCEVdj:thm:strong_convSDDE}
Suppose Assumptions A and B hold. Then the delay semi-discrete method $y_t^0,$ with $m=1/4$ converges in the mean square sense to the true solution of the delay CIR/CEV model, i.e. (\ref{CIRCEVdj:eq:sddej})  with $g\equiv0,$ with order of convergence given by

\beqq \label{CIRCEVdj:eq:strong_convSDDE_pol}
\bfE\sup_{0\leq t\leq T}|y_t^0-x_t^0|^2\leq \left\{\barr{ll}
C\D^{(\alpha-1/2)\wedge\gamma} & , \alpha\in(1/2,1),\\
&\\
C\D^{\gamma\wedge(1/2)}& ,  \alpha=1/2,
\earr \right.
\eeqq
where $C$ is a positive constant independent of $\D.$
\ethe

The (JASDM) method reads
\beao
y_{t_{k+1}} &=& y_{t_{{k+1}^-}} + g(y_{t_{{k+1}^-}}) \D\wt{N}_k\\
&=& y_{t_{{k+1}^-}} + g(y_{t_{{k+1}^-}}) (\bbi_{\{ t_{k+1} \hbox{ is a jump time}\}}-\lam\D_k),
\eeao
where $y_{t_{{k+1}^-}}$ approximates $x_{t_{{k+1}^-}}.$ Assumption A implies for any $r>0$, 
$$r+g(r)(1-\lam\D_k)\geq r-(1-\lam\D_k)Lr>r(1-L)>0,$$
when $L<1$ and $r-g(r)\lam\D_k\geq r(1-\lam L\D_k)>0,$ for any $\D<(\lam L)^{-1}.$ In case $g(\cdot)>0$ 
$$r+g(r)(1-\lam\D_k)>g(r)(1-\lam\D_k)>0,$$ for any $\D<\lam^{-1}.$ Therefore, $y_{t_{k+1}}>0$ a.s. when $\D<\lam^{-1}(L^{-1}\wedge1)$.

The results of Theorem \ref{CIRCEVdj:thm:strong_convSDDE} are valid in the case jumps are included in the model and this is stated in the following result.

\bth[Strong convergence of (JASDM) to (\ref{CIRCEVdj:eq:sddej})] \label{CIRCEVdj:thm:strong_convSDDEjump}
Suppose Assumption A and B hold and $\D<\lam^{-1}(L^{-1}\wedge1)$. Then the jump adapted semi-discrete method (\ref{CIRCEVdj:eq:JASD}) converges in the mean square sense to the true solution of the jump-extended CIR/CEV model (\ref{CIRCEVdj:eq:sddej}) with order of convergence given by

\beqq \label{CIRCEVdj:eq:strong_convSDDEjump_pol}
\bfE\sup_{0\leq t\leq T}|y_t^g-x_t^g|^2\leq \left\{\barr{ll}
C\D^{(\alpha-1/2)\wedge\gamma} & , \alpha\in(1/2,1),\\
&\\
C\D^{\gamma\wedge(1/2)}& ,  \alpha=1/2,
\earr \right.
\eeqq
where $C$ is a positive constant independent of $\D.$
\ethe

\section{Numerical illustration}\label{CIRCEVdj:sec:numerics}

\begin{figure}[ht]
	\centering
	\begin{subfigure}{.45\textwidth}
		\includegraphics[width=1\textwidth]{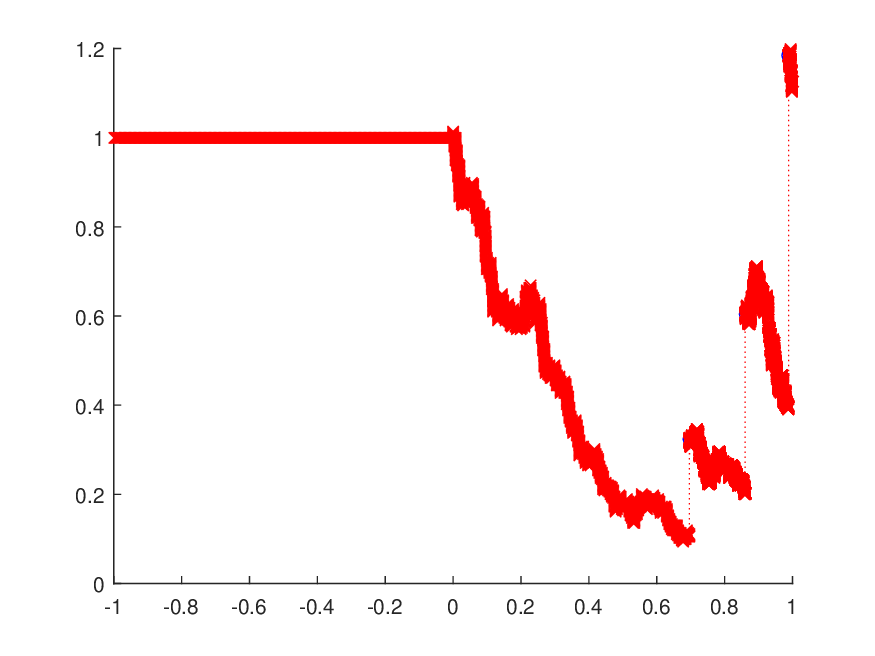}\label{CIRCEVdj-fig:exp1ex1}
	\end{subfigure}
	\begin{subfigure}{.45\textwidth}
		\includegraphics[width=1\textwidth]{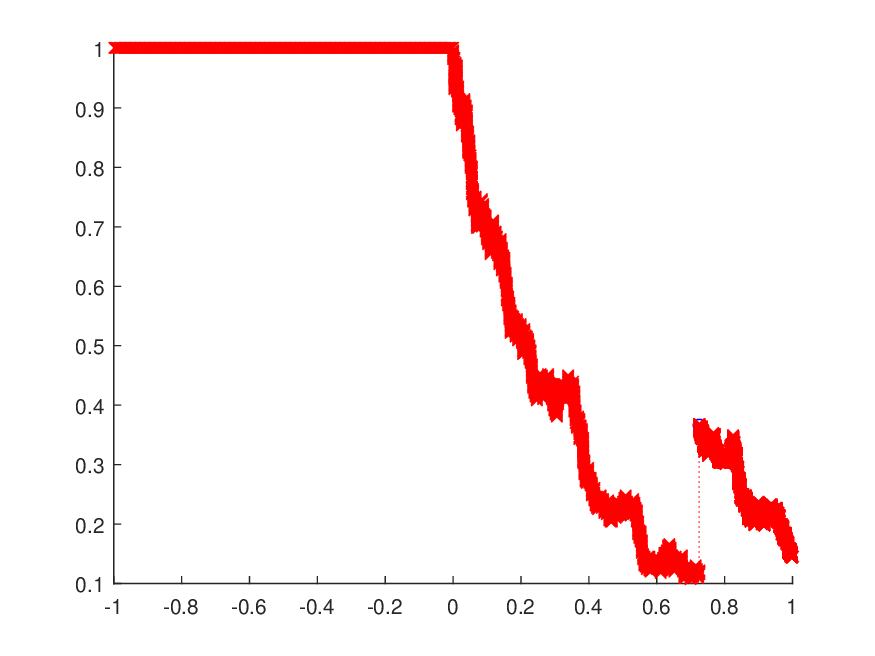}\label{CIRCEVdj-fig:exp1ex2}
	\end{subfigure}
	\caption{Two one-path simulations of solution of (\ref{CIRCEVdj:eq:sddejEX}), SETI with $\alpha$=0.5, $\gamma$=1, $\D$=$2^{-14}$.}\label{CIRCEVdj-fig:exp1}
\end{figure}

\begin{figure}[ht]
	\centering
	\begin{subfigure}{.45\textwidth}
		\includegraphics[width=1\textwidth]{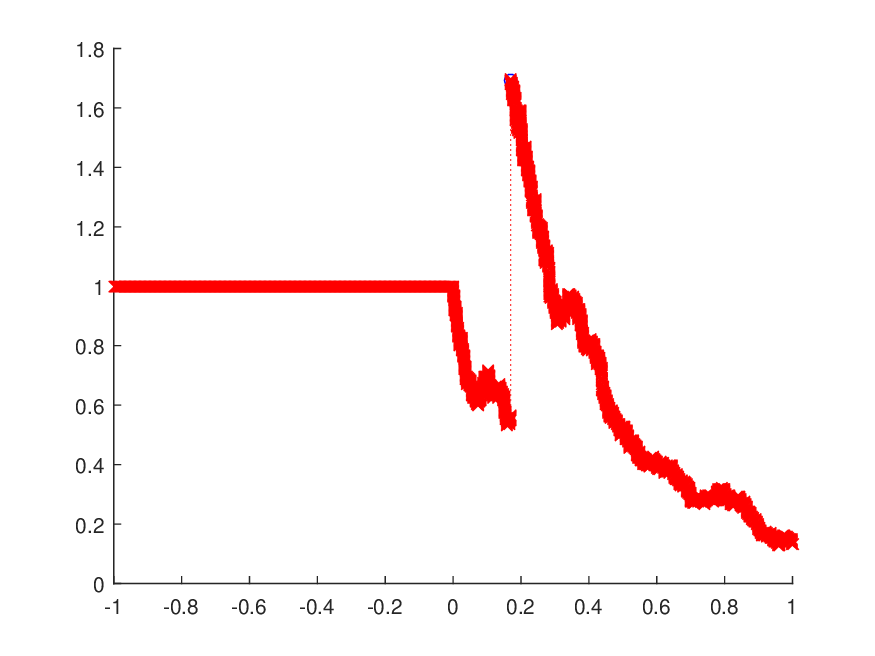}\label{CIRCEVdj-fig:exp2ex1}
	\end{subfigure}
	\begin{subfigure}{.45\textwidth}
		\includegraphics[width=1\textwidth]{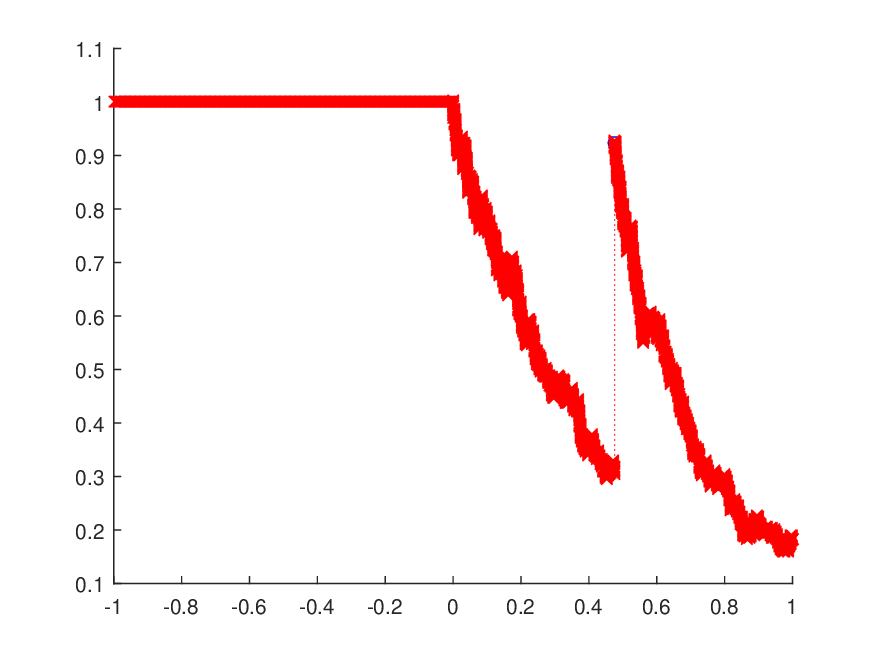}\label{CIRCEVdj-fig:exp2ex2}
	\end{subfigure}
	\caption{Two one-path simulations of solution of (\ref{CIRCEVdj:eq:sddejEX}), SETI with $\alpha$=0.7, $\gamma$=1/2, $\D$=$2^{-14}$.}\label{CIRCEVdj-fig:exp2}
\end{figure}

\begin{figure}[ht]
	\centering
	\begin{subfigure}{.45\textwidth}
		\includegraphics[width=1\textwidth]{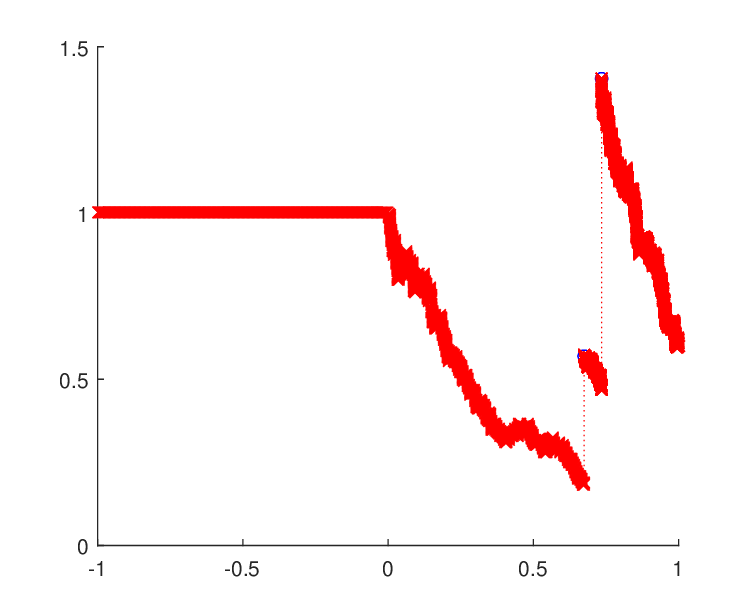}\label{CIRCEVdj-fig:exp3ex1}
	\end{subfigure}
	\begin{subfigure}{.45\textwidth}
		\includegraphics[width=1\textwidth]{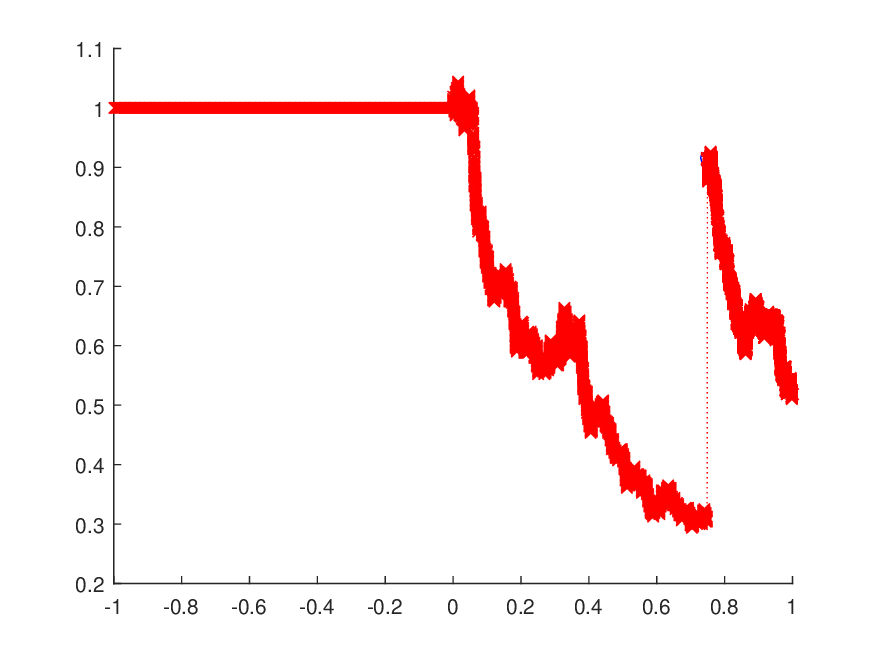}\label{CIRCEVdj-fig:exp3ex2}
	\end{subfigure}
	\caption{Two one-path simulations of solution of (\ref{CIRCEVdj:eq:sddejEX}), SETI with $\alpha$=0.9, $\gamma$=1/2, $\D$=$2^{-14}$.}\label{CIRCEVdj-fig:exp3}
\end{figure}

\begin{figure}[ht]
	\centering
	\begin{subfigure}{.45\textwidth}
		\includegraphics[width=1\textwidth]{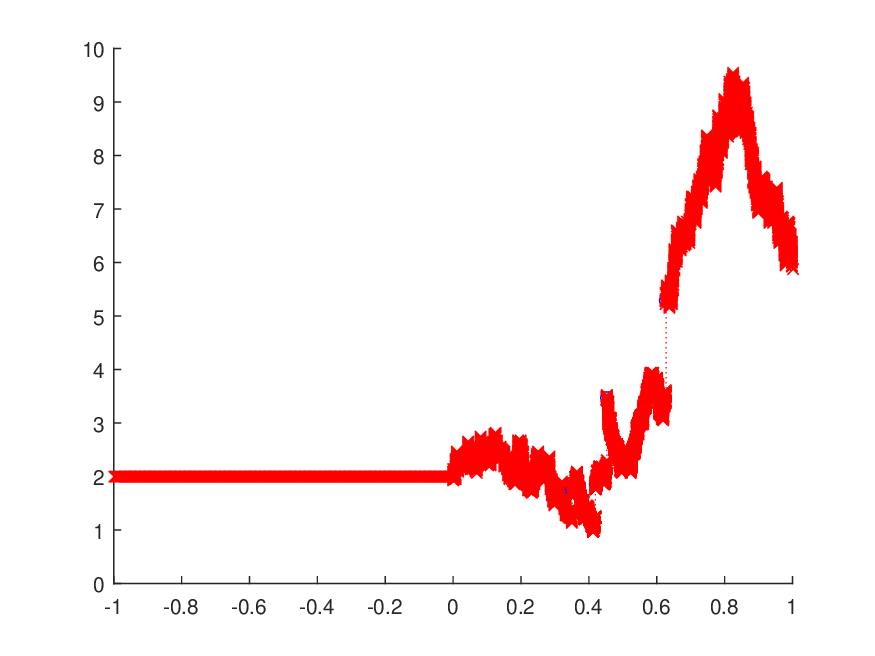}\label{CIRCEVdj-fig:SETIIex1}
	\end{subfigure}
	\begin{subfigure}{.45\textwidth}
		\centering
		\includegraphics[width=1\textwidth]{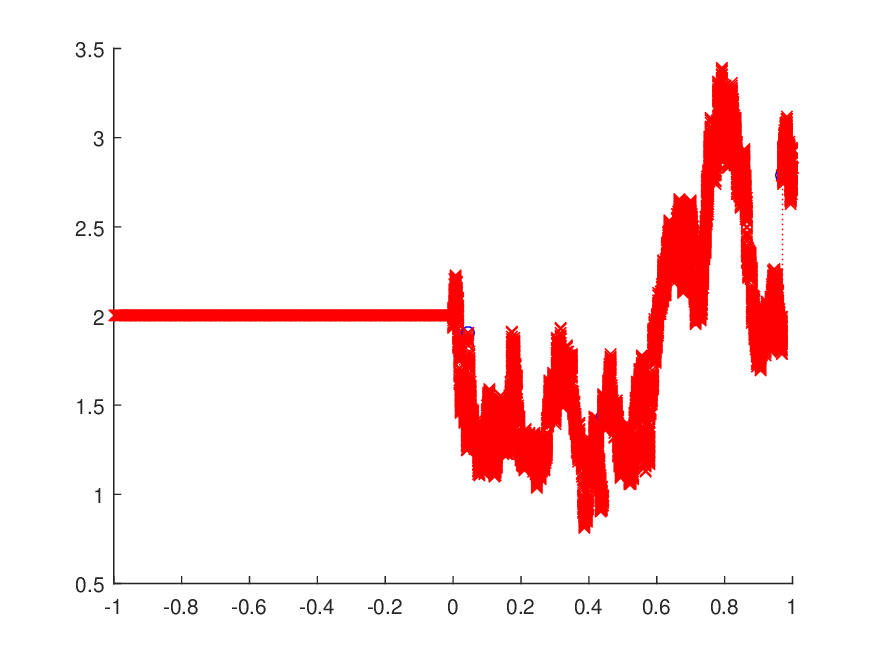}\label{CIRCEVdj-fig:SETIIex2}
	\end{subfigure}
	\caption{Two one-path simulations of solution of (\ref{CIRCEVdj:eq:sddejEX}), SETII with  $\D$=$2^{-14}$.}\label{CIRCEVdj-fig:SETII}
\end{figure}

In this section, we present some numerical results that illustrate the positivity preserving and the order of strong convergence of (JASDM). We examine models with coefficients as in \cite[Sec. ~4]{yang_etal:2018} where CIR/CEV jump models are considered with no delay. We 
consider the following SDDEs with jumps,
\beam
\label{CIRCEVdj:eq:sddejEX}
dx_t &=& (k_1 - k_2x_{t^-})dt + k_3b(x_{t-1}) (x_{t^-})^\alpha dW_t + g(x_{t^-})d\wt{N}_t,  \quad t\in [0,1],\\
\nonumber x_t& = &\xi(t),\quad  t\in[-1,0],
\eeam
that is model (\ref{CIRCEVdj:eq:sddej}) with parameters from Table \ref{CIRCEVdj:table:t2}.
\begin{table}[h!]
	\begin{center}
		\caption{Parameters in model (\ref{CIRCEVdj:eq:sddej}).}
		\label{CIRCEVdj:table:t2}
		\begin{tabular}{r|r|r|r|r|r|r|r|r|r|r|r} 
			& $k_1$  & $k_2$ & $k_3$ & $\alpha$ & $b(x)$ & $\gamma$ & $g(x)$ & $\xi(t)$ & $\lam$ & $\tau$ & T\\
			\hline
			SETI & $0.24$  & $3$ & $0.4$ & $0.5, 0.7, 0.9$ & $x^\gamma$ & $1/2, 1$ & $2x$ & $1$ &1 & 1& 1\\
			\hline
			SETII &$2$  &$2$ & $1.5$ & $0.5$ & $1+e^{-x}$ & $1$ & $0.5x$ & $2$ & 1 & 1& 1\\
		\end{tabular}
	\end{center}
\end{table}   
Let  $\bbt$ be a jump-adapted time partition of the interval $[-1,1].$ The (JASDM) method for $\theta=1/2$ reads $y_t=\xi(t)$ when $t\in[-1,0]$ and for $k=0, 1, \ldots, n_T-1,$ see (\ref{CIRCEVdj:eq:JASD})
\beqq\label{CIRCEVdj:eq:JASDeX}
\left\{\barr{ll}
y_{t_{{k+1}^-}} = (z_{t_{k+1}})^2,\\
y_{t_{k+1}} = y_{t_{{k+1}^-}} +g(y_{t_{{k+1}^-}})(\D N_k-\D_k),\\
\earr \right.
\eeqq
where 
\beao
z_{t} &=&  \sqrt{y_{t_k}\left(1- \frac{2k_2\D_k}{2+k_2\D_k}\right) + \frac{2k_1\D_k}{2+k_2\D_k} - \frac{(k_3)^2}{(2+k_2\D_k)^2 }\frac{b^2(y_{t_{k-1}})}{(1+b(y_{t_{k-1}})\D_k^{1/4})^2}(y_{t_k})^{2\alpha-1}\D_k}\\
&& +\frac{k_3}{1+k_2\D_k} \frac{b(y_{t_{k-1}})}{1+b(y_{t_{k-1}})\D_k^{1/4}}(y_{t_k})^{\alpha -\frac{1}{2}}(W_t - W_{t_k})
\eeao
and $\D_k=t_{k+1}-t_{k}.$ First, we plot two one-path simulations of (\ref{CIRCEVdj:eq:sddejEX}) for different parameters of Table \ref{CIRCEVdj:table:t2}. We use (JASDM) with $\D=2^{-14}$ as the exact solution, see Figures \ref{CIRCEVdj-fig:exp1}, \ref{CIRCEVdj-fig:exp2},  \ref{CIRCEVdj-fig:exp3} and \ref{CIRCEVdj-fig:SETII}.

\begin{figure}[ht]
	\centering
	\begin{subfigure}{.45\textwidth}
		\includegraphics[width=1\textwidth]{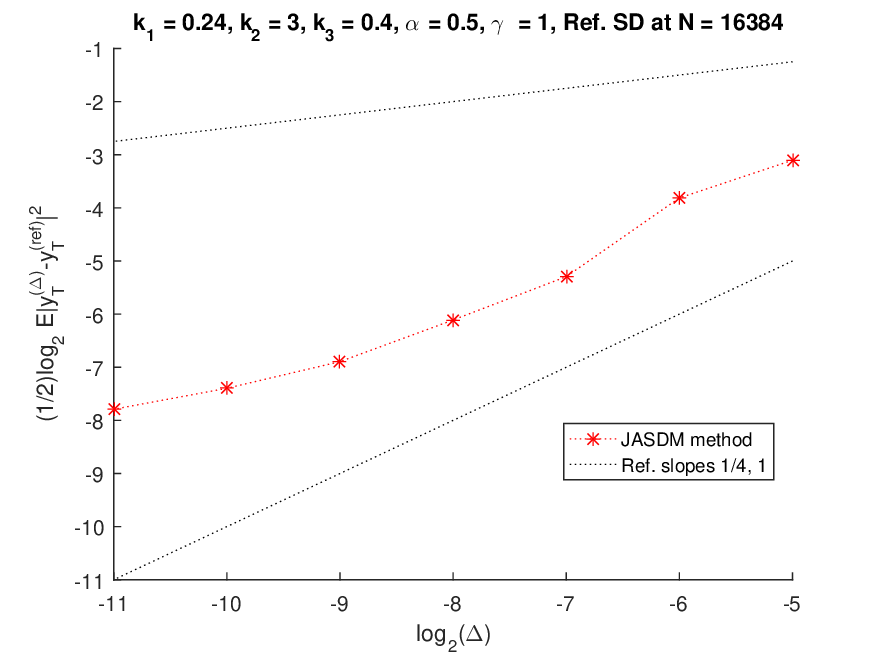}\label{CIRCEVdj-fig:CIRdelayJ1_REP}
		\caption{SET I}
	\end{subfigure}
	\begin{subfigure}{.45\textwidth}
		\includegraphics[width=1\textwidth]{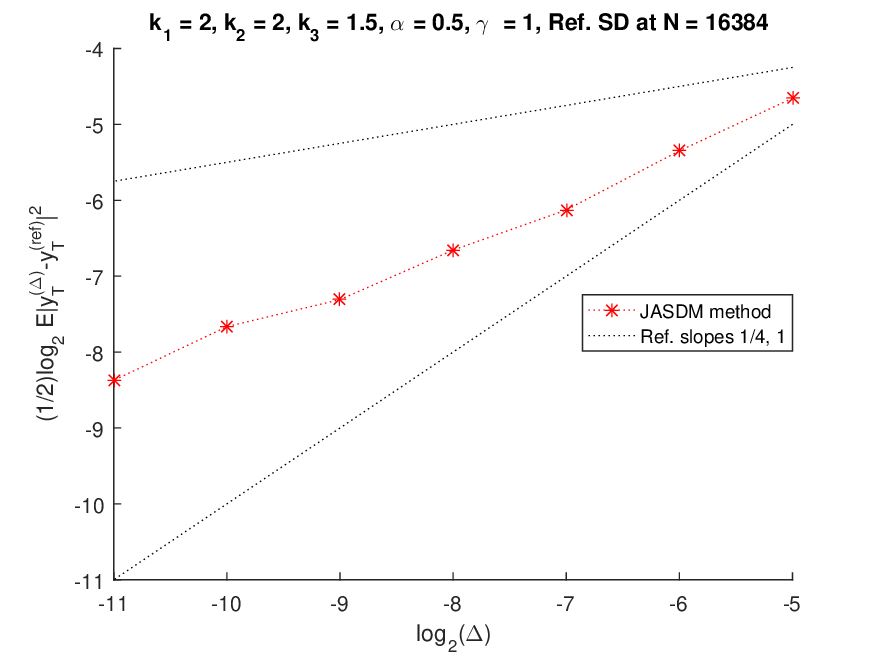}\label{CIRCEVdj-fig:CIRdelayJ_SETII}
		\caption{SET II}
	\end{subfigure}
	\caption{Error and step size for (JASDM) for CIR model (\ref{CIRCEVdj:eq:sddejEX}) with $\gamma=1.$}\label{CIRCEVdj-fig:CIRdelayJ12}
\end{figure}

\begin{figure}[ht]
	\centering
	\begin{subfigure}{.45\textwidth}
		\includegraphics[width=1\textwidth]{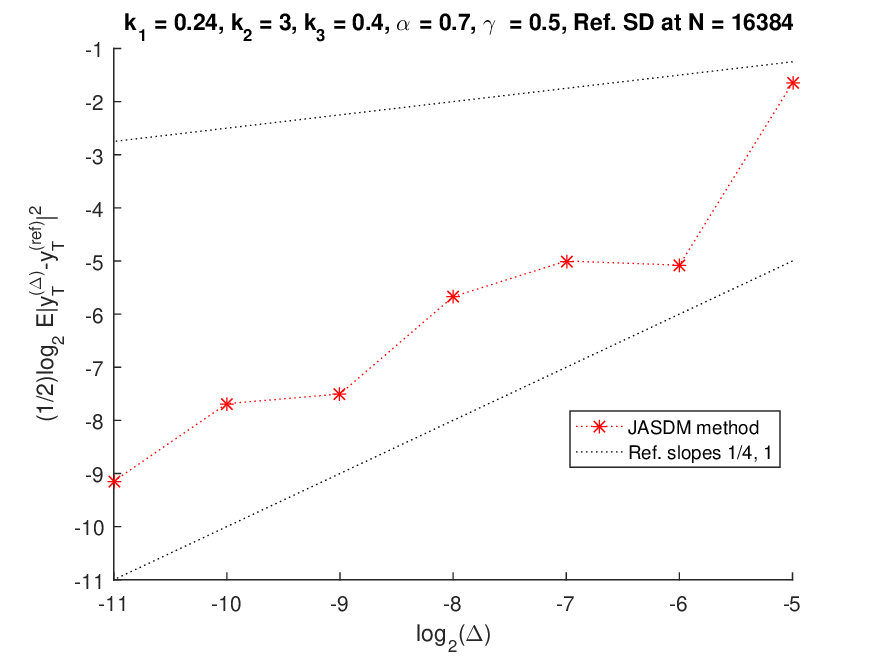}\label{CIRCEVdj-fig:CEVdelayJ1}
		\caption{Case $(\alpha,\gamma)=(0.7,1/2).$}
	\end{subfigure}
	\begin{subfigure}{.45\textwidth}
		\includegraphics[width=1\textwidth]{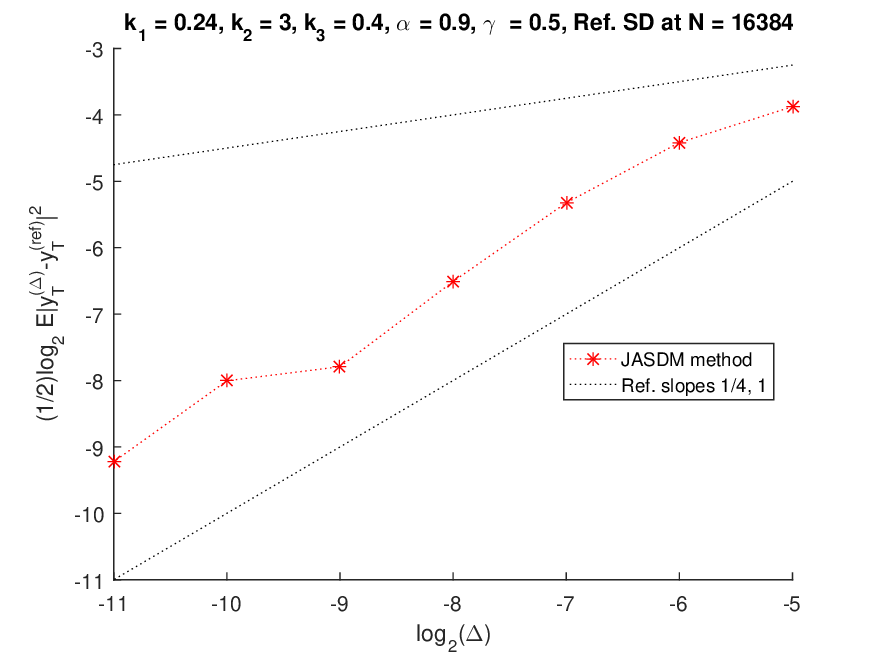}\label{CIRCEVdj-fig:CEVdelayJ2}
		\caption{Case $(\alpha,\gamma)=(0.9,1/2).$}
	\end{subfigure}
	\caption{Error and step size for (JASDM) for CEV model (\ref{CIRCEVdj:eq:sddejEX}) for different $\alpha.$}\label{CIRCEVdj-fig:CEVdelayJ12}
\end{figure}

We estimate the endpoint $\bbl^2$-norm $\ep=\sqrt{\bfE|y^{(\D)}(T) - x_T|^2},$ of the difference between the (JASDM) method at step size $\D$ and the exact solution of (\ref{CIRCEVdj:eq:sddejEX}). We compute $M$ batches of $L$ simulation paths, where each batch is estimated by 
$\hat{\ep}_j=\frac{1}{L}\sum_{i=1}^L|y_{i,j}^{(\D)}(T) - y_{i,j}^{(ref)}(T)|^2$ 
and the Monte Carlo estimator of the error is
\beqq\label{CIRCEVdj:l2error}
\hat{\ep}=\sqrt{\frac{1}{ML}\sum_{j=1}^M\sum_{i=1}^L|y_{i,j}^{(\D)}(T) - y_{i,j}^{(ref)}(T)|^2}
\eeqq
and requires $M\cdot L$ Monte Carlo sample paths.  We also use $M\cdot L$ Poisson paths. The reference solution is calculated using the method at a very fine time grid, $\D=2^{-14}.$ We have shown in Theorem \ref{CIRCEVdj:thm:strong_convSDDEjump} that the (JASDM) converges strongly to the exact solution, so we use the (JASDM) method as a reference solution. 
We simulate $50\cdot 100=5000$ paths, where the choice of the number of Monte Carlo paths  is adequately large, so as not to significantly hinder the mean-square errors. We compute the approximation error (\ref{CIRCEVdj:l2error}) and present the results in Table \ref{CIRCEVdj:table:t3}. We also give a graphical illustration for the delay CIR jump models in a $\log_2$-$\log_2$ scale in Figure \ref{CIRCEVdj-fig:CIRdelayJ12}  and for the delay CEV jump models in Figure \ref{CIRCEVdj-fig:CEVdelayJ12}.
\footnotesize
\begin{table}[h!]
	\begin{center}
		\caption{Convergence results for (JASDM), see (\ref{CIRCEVdj:eq:JASDeX}).}
		\label{CIRCEVdj:table:t3}
		\begin{tabular}{r|rr|rr|rr|rr} 
			& $\alpha=1/2$ & SETI    & $\alpha=0.7$ & SETI  &  $\alpha=0.9$ & SETI  & $\alpha=1/2$ & SETII\\
			$\D$ & $\gamma=1$ & rate   & $\gamma=1/2$ & rate &  $\gamma=1/2$ & rate  &$\gamma=1$ & rate \\
			\hline
			$2^{-5}$ & $0.1163$ & $-$                    & $0.3188$ &   $-$               & $0.0681$     & $-$              & $0.0397$ 	 	 & $-$  \\
			$2^{-6}$ & $0.0712$ & $0.7084$   & $0.0296$ & $3.4287$ & $0.0468$ & $0.5423$ & $0.0245$ & $0.6931$ \\
			$2^{-7}$ & $0.0255$ & $1.4824$   & $0.0312$ & $-0.0777$ & $0.0250$ & $0.9053$ & $0.0143$ & $0.7792$\\
			$2^{-8}$ & $0.0145$ & $0.8157$   & $0.0196$ & $0.6716$ & $0.0110$ & $1.1851$ & $0.0099$ & $0.5287$\\
			$2^{-9}$ & $0.0084$ & $0.7838$   & $0.0055$ & $1.8319$ & $0.0045$ & $1.2834$& $0.0063$ & $0.6524$\\
			$2^{-10}$ &$0.0059$ & $0.5019$  & $0.0049$ & $0.1818$ & $0.0039$ & $0.2120$ & $0.0049$ & $0.3564$\\
			$2^{-11}$ &$0.0045$ & $0.3944$    & $0.0018$ & $1.4601$ & $0.0017$ & $1.2191$ & $0.0030$ & $0.7078$
		\end{tabular} 
	\end{center}
\end{table}\normalsize   

Assumption B requires $\D<0.42\wedge0.18\wedge0.64=0.18$ therefore step $\D=2^{-5}$ is sufficient. Moreover, $\lam=L=1$ therefore the other requirement for the step size reads $\D<1,$ see Theorem \ref{CIRCEVdj:thm:strong_convSDDEjump}. For the delay CIR jump model with $\gamma =1,$ the convergence rate should be at least $1/4$ and for the CEV type models with $(\alpha,\gamma)=(0.7,1/2)$ and $(\alpha,\gamma)=(0.9,1/2)$ the convergence rate should be at least $0.1$ and $0.2$ respectively.

\section{Proofs}\label{CIRCEVdj:sec:proofs}

In this section we provide the proofs of all the results of Section \ref{CIRCEVdj:sec:main}. 
\bpf[Proof of Lemma \ref{CIRCEVdj:lem:moment_bound}]
We define for any $k>0$ the stopping time
$$
\tau_k:=\inf\{t\in [0,T]: |x_t|>k\}
$$
for the stochastic process $x_t$  of (\ref{CIRCEVdj:eq:sddejABS}), with the convection that $\inf\emptyset=\infty.$ The H\"{o}lder inequality and the Burkholder-Davis-Gundy inequality, see \cite{platen:2010}, imply for any $t_1\in[0,T]$ and $p>2,$
\beam\nonumber
\bfE\sup_{0\leq t\leq t_1}|x_{t\wedge \tau_k}|^p&\leq& 4^{p-1}\bfE |\xi_0|^p + 4^{p-1}\bfE\left(\int_{0}^{t_1\wedge \tau_k}|k_1 - k_2x_{s^-}|ds\right)^p\\
&&\nonumber+ 4^{p-1}(k_3)^p\bfE\left(\sup_{0\leq t\leq t_1}\left|\int_{0}^{t\wedge \tau_k}b(|x_{s-\tau}|)|x_{s^-}|^\alpha dW_s\right|^p\right)\\
&&\nonumber + 4^{p-1}\bfE\left(\sup_{0\leq t\leq t_1}\left|\int_{0}^{t\wedge \tau_k}g(x_{s^-}) d\wt{N}_s\right|^p\right)\\
&\leq&\nonumber 4^{p-1}\bfE |\xi_0|^p + 4^{p-1}T^{p-1}\bfE\int_{0}^{t_1\wedge \tau_k}|k_1 - k_2x_{s^-}|^pds\\
&&\nonumber+ 4^{p-1}(k_3)^pC_p\bfE\left(\int_{0}^{t_1\wedge \tau_k}b(|x_{s-\tau}|)^2|x_{s^-}|^{2\alpha} ds\right)^{p/2}\\
&&\label{CIRCEVdj:eq:mb1} + 4^{p-1}C_{p,\lam}\bfE\left(\int_{0}^{t_1\wedge \tau_k}|g(x_{s^-})|^p ds\right),
\eeam
where $C_p, C_{p,\lam}$ are constants depending on $p$ and $p,\lam$ respectively. The first integral in (\ref{CIRCEVdj:eq:mb1}) is bounded in the following way
\beqq\label{CIRCEVdj:eq:mb2}
\bfE\int_{0}^{t_1\wedge \tau_k}|k_1 - k_2x_{s^-}|^pds\leq 2^{p-1}(k_1)^pT^p + 
2^{p-1}(k_2)^p\int_{0}^{t_1}\bfE| x_{s^-\wedge \tau_k}|^pds.
\eeqq
The third integral in (\ref{CIRCEVdj:eq:mb1}) can be estimated using Assumption A, as
\beqq\label{CIRCEVdj:eq:mb3}
\bfE\left(\int_{0}^{t_1\wedge \tau_k}|g(x_{s^-})|^p ds\right)\leq 2^{p-1}L^p\bfE\int_{0}^{t_1}\bfE|(x_{s^-\wedge \tau_k})|^p ds + 2^{p-1}|g(0)|^pT.
\eeqq
Moreover, by the elementary inequality $u^{1-\beta}v^{\beta}\leq(1-\beta)u + \beta v,$ for $u,v\geq0$ and $0\leq\beta<1$ and another application of H\"{o}lder's inequality we get\footnotesize
\beao
&&\bfE\left(\int_{0}^{t_1\wedge \tau_k}b(|x_{s-\tau}|)^2|x_{s^-}|^{2\alpha} ds\right)^{p/2}\leq \bfE\left(\int_{0}^{t_1\wedge \tau_k}(1-\alpha)b(|x_{s-\tau}|)^{2/(1-\alpha)} + \alpha|x_{s^-}|^{2} ds\right)^{p/2}\\
&\leq& 2^{\frac{p}{2}-1}\left[\bfE\left(\int_{0}^{t_1\wedge \tau_k}(1-\alpha)b(|x_{s-\tau}|)^{2/(1-\alpha)}ds\right)^{p/2} + \bfE\left(\int_{0}^{t_1\wedge \tau_k}\alpha|x_{s^-}|^{2} ds\right)^{p/2}\right]\\
&\leq& 2^{\frac{p}{2}-1}(1-\alpha)^{p/2}T^{\frac{p}{2}-1}\int_{0}^{t_1}\bfE b(|x_{s\wedge \tau_k-\tau}|)^{p/(1-\alpha)}ds +   2^{\frac{p}{2}-1}\alpha^{p/2}T^{\frac{p}{2}-1}\int_{0}^{t_1}\bfE|x_{s^-\wedge \tau_k}|^{p} ds\\
&\leq& C + (2T(1-\alpha))^{\frac{p}{2}-1}C_\gamma^{\frac{p}{1-\alpha}}\int_{0}^{t_1}(1-\alpha)\bfE|x_{s\wedge \tau_k-\tau}|^{\frac{\gamma}{1-\alpha}p}ds +   (2T\alpha)^{\frac{p}{2}-1}\int_{0}^{t_1}\alpha\bfE|x_{s^-\wedge \tau_k}|^{p} ds,
\eeao\normalsize
where the constant $C$ depends on $\bfE (|\xi_{-\tau}|^{\frac{p}{1-\alpha}})$ and also in the last step we used Assumption A. Plugging the last estimate and (\ref{CIRCEVdj:eq:mb2}), (\ref{CIRCEVdj:eq:mb3}) into (\ref{CIRCEVdj:eq:mb1}) yields

\beqq\label{CIRCEVdj:eq:mb4}
\bfE\sup_{0\leq t\leq t_1}|x_{t\wedge \tau_k}|^p\leq \bba_p + \bbb_p\int_{0}^{t_1}\bfE|x_{s\wedge \tau_k-\tau}|^{\frac{\gamma}{1-\alpha}p}ds + \bbc_p\int_{0}^{t_1}\bfE|x_{s^-\wedge \tau_k}|^{p}ds,
\eeqq
where the constants $\bba_p,\bbb_p$ and $\bbc_p$ depend on $p, \bfE\|\xi\|, k_1, k_2, k_3, T,\lam$ and $\alpha.$ 

We consider a finite sequence $\{p_1, p_2, \ldots, p_{[T/\tau] + 1}\}$ such that $p_i>2$ and $\gamma p_{i+1}<(1-\alpha)p_i$ for $i =1, 2,\ldots, [T/\tau]$ where $[z]$ denotes the integer part of $z.$ We work successively in the regions $[0,\tau], [0,2\tau],$ e.t.c. until we cover the interval $[0,T]$ proving in each case the moment bounds.   This idea can be found in the proofs of \cite[Lemma 2.1]{jiang_etal:2011} or \cite[Theorem 2.1]{wu_etal:2009}. Take a $t_1\in[0,\tau].$ Then $\bfE |x_{s-\tau}|^{\frac{\gamma}{1-\alpha}p_1}\leq\bfE \| \xi\|^{\frac{\gamma}{1-\alpha}p_1}$  for $s\in[0,t_1].$ Equation (\ref{CIRCEVdj:eq:mb4}) yields
\beao
\bfE\sup_{0\leq t\leq t_1}|x_{t\wedge \tau_k}|^{p_1}&\leq& \bba_{p_1} + \bbb_{p_1}T\bfE\| \xi \|^{\frac{\gamma}{1-\alpha}p_1} + \bbc_{p_1}\int_{0}^{t_1}\bfE\sup_{0\leq r\leq s}|x_{r\wedge \tau_k}|^{p_1}ds\\
&\leq& (\bba_{p_1} + \bbb_{p_1}T\bfE\| \xi \|^{\frac{\gamma}{1-\alpha}p_1})e^{\bbc_{p_1}\tau},
\eeao
by application of the Gronwall inequality, which in turn implies
\beqq
\bfE\sup_{0\leq t\leq t_1}|x_{t}|^{p_1}\leq (\bba_{p_1} + \bbb_{p_1}T\bfE\| \xi \|^{\frac{\gamma}{1-\alpha}p_1})e^{\bbc_{p_1}\tau},
\eeqq
by application of Fatou's theorem as we let $\kto.$
We have 
\beao
\bfE\sup_{-\tau\leq t\leq \tau}|x_{t}|^{p_1}&\leq& \bfE\sup_{-\tau\leq t\leq 0}|x_{t}|^{p_1} + \bfE\sup_{0\leq t\leq \tau}|x_{t}|^{p_1}\\
&\leq& \bfE\| \xi \|^{p_1} +   (\bba_{p_1} + \bbb_{p_1}T\bfE\| \xi \|^{\frac{\gamma}{1-\alpha}p_1})e^{\bbc_{p_1}\tau}:=A_1.
\eeao
When  $t_1\in[0,2\tau]$ we get by (\ref{CIRCEVdj:eq:mb4}) 
\beao
\bfE\sup_{0\leq t\leq t_1}|x_{t\wedge \tau_k}|^{p_2}&\leq& \bba_{p_2} + \bbb_{p_2}
\int_0^{2\tau}\left(\bfE|x_{s\wedge \tau_k-\tau}|^{p_1}\right)^{\frac{\gamma p_2}{(1-\alpha)p_1}}ds + \bbc_{p_2}\int_{0}^{t_1}\bfE\sup_{0\leq r\leq s}|x_{r^-\wedge \tau_k}|^{p_2}ds\\
&\leq& \bba_{p_2} + \bbb_{p_2}
\int_{-\tau}^{\tau}\left(\bfE|x_{s\wedge \tau_k}|^{p_1}\right)^{\frac{\gamma p_2}{(1-\alpha)p_1}}ds  + \bbc_{p_2}\int_{0}^{t_1}\bfE\sup_{0\leq r\leq s}|x_{r^-\wedge \tau_k}|^{p_2}ds\\
&\leq&\bba_{p_2} + 2\bbb_{p_2}\tau A_1^{\frac{\gamma p_2}{(1-\alpha)p_1}} + \bbc_{p_2}\int_{0}^{t_1}\bfE\sup_{0\leq r\leq s}|x_{r^-\wedge \tau_k}|^{p_2}ds\\
&\leq&(\bba_{p_2} + 2\bbb_{p_2}TA_1^{\frac{\gamma p_2}{(1-\alpha)p_1}})e^{2\bbc_{p_2}\tau}
\eeao
and once again by the Fatou theorem we have
$$
\bfE\sup_{0\leq t\leq 2\tau}|x_{t}|^{p_2}\leq(\bba_{p_2} + 2\bbb_{p_2}TA_1^{\frac{\gamma p_2}{(1-\alpha)p_1}})e^{2\bbc_{p_2}\tau}.
$$
and analogously 
\beao
\bfE\sup_{-\tau\leq t\leq 2\tau}|x_{t}|^{p_2}&\leq& \bfE\sup_{-\tau\leq t\leq 0}|x_{t}|^{p_2} + \bfE\sup_{0\leq t\leq 2\tau}|x_{t}|^{p_2}\\
&\leq& \bfE\| \xi \|^{p_2} +   (\bba_{p_2} + 2\bbb_{p_2}TA_1^{\frac{\gamma p_2}{(1-\alpha)p_1}})e^{2\bbc_{p_2}\tau}:=A_2.
\eeao
Repeating the previous steps we finally obtain that for any $p=p_{[T/\tau]+1}>2$ there exists a positive constant $A_{[T/\tau]+1}$ such that 
$$
\bfE\sup_{-\tau\leq t\leq T}|x_{t}|^{p}\leq A_{[T/\tau]+1}:=A.
$$
The case $0<p\leq2$ follows from the H\"{o}lder inequality since then $\bfE |x_t|^p\leq (\bfE |x_t|^3)^{p/3}\leq A^{p/3}.$ Therefore (\ref{CIRCEVdj:eq:moment_bound}) is true.
\epf

\bpf[Proof of Proposition \ref{CIRCEVdj:prop:uniq_posABS}]
The uniqueness of the solution of (\ref{CIRCEVdj:eq:sddejABS}) can be found in \cite{maoX:1994}. The non negativity property can be shown by adopting the Yamada-Watanabe approach.
Let the non-increasing sequence $\{e_m\}_{m\in\bbN}$ with $e_m=e^{-m(m+1)/2}$ and $e_0=1.$ We introduce the following sequence $\phi_m(x)=0$ for $x\geq0$ and
$$
\phi_m(x)=\int_0^{-x}dy\int_0^{y}\psi_m(u)du, \quad \hbox{ for } \, x<0,
$$
where the existence of the continuous function $\psi_m(u)$ with $0\leq \psi_m(u) \leq 2/(mu)$ and support in $(e_m,e_{m-1})$ is justified by $\int_{e_m}^{e_{m-1}}(du/u)=m.$ The following relations hold for $\phi_m\in\bbc^2(\bbR,\bbR)$ with $\phi_m(0)=0,$
$$
x^- - e_{m-1}\leq\phi_m(x)\leq x^-, \hbox{ for  } \quad x\in\bbR, 
$$
where $x^-=-x\vee0,$
$$
-1\leq\phi_{m}^{\prime}(x)\leq0, \,\hbox{ if }  x<-e_m \,\hbox{ and }  \,  \phi_{m}^{\prime}(x)=0 \,\hbox{ otherwise;} $$
$$
|\phi_{m}^{\prime \prime }(x)|\leq\frac{2}{m|x|}, \,\hbox{ when }  \,-e_{m-1}<x<-e_{m} \,\hbox{ and }  \,  \phi_{m}^{\prime \prime }(x)=0 \,\hbox{ otherwise.}
$$
Application of the It\^o formula for the jump-extended CIR/CEV model reads\footnotesize
\beao
\phi_{m}(x_t) &=& \phi_{m}(x_0) + \int_0^t \phi_{m}^{\prime}(x_s)(k_1 - k_2x_{s^-}-\lam g(x_{s^-}))ds + \frac{(k_3)^2}{2}\int_{0}^{t} \phi_{m}^{\prime\prime}(x_s)b(|x_{s-\tau}|)^2|x_{s^-}|^{2\alpha}ds\\
&& + k_3\int_{0}^{t} \phi_{m}^{\prime}(x_s)b(|x_{s-\tau}|)|x_{s^-}|^{\alpha}dW_s + \int_0^t \left(\phi_{m}(x_{s^-} +g(x_{s^-})) - \phi_{m}(x_{s^-}) \right)dN_s.
\eeao\normalsize
Now, we take expectations in the above inequality, using the moment bounds of Lemma \ref{CIRCEVdj:lem:moment_bound} and the properties of $\phi_m$ to get in the case $g(0)=0$ 
\footnotesize
\beao
\bfE\phi_{m}(x_t) &=&\int_0^t \bfE\phi_{m}^{\prime}(x_s)(k_1 - k_2x_{s^-}-\lam g(x_{s^-}))ds + \frac{(k_3)^2}{2}\bfE\left(\int_{0}^{t} \phi_{m}^{\prime\prime}(x_s)b(|x_{s-\tau}|)^2|x_{s^-}|^{2\alpha}ds\right)\\
&& + \lam\bfE\int_0^t \left(\phi_{m}(x_{s^-} +g(x_{s^-})) - \phi_{m}(x_{s^-}) \right)ds\\
&\leq& \frac{(k_3C_\gamma)^2}{m}\int_{0}^{t}\bfE \left(|x_{s-\tau}|^{2\gamma}|x_{s^-}|^{2\alpha-1}\right)ds + \frac{(k_3)^2}{m}\int_{0}^{t}\bfE \left(b|(\xi_{-\tau})^2||x_{s^-}|^{2\alpha-1}\right)ds\\
&& + k_2\sup_{x\in\bbR}|\phi_{m}^{\prime}(x)|\int_0^t \bfE (x_{s^-})ds + \lam\sup_{x\in\bbR}|\phi_{m}^{\prime}(x)|\int_0^t \bfE( g(x_{s^-}))ds\\ 
&\leq& \frac{C(k_3)^2}{m}\int_{0}^{t}\sqrt{\bfE |x_{s-\tau}|^{4\gamma}}
\sqrt{\bfE |x_{s^-}|^{4\alpha-2}}ds + (\lam L + k_2)\int_0^t \bfE (x_{s^-}) ds\\ 
&\leq& \frac{C(k_3)^2}{m}A + (\lam L + k_2)\int_0^t (\bfE (\phi_m(x_{s^-})) + e_{m-1})ds,
\eeao\normalsize
where $C$ is a constant or 
\beao
\bfE\phi_{m}(x_t) &\leq&  \left(\frac{C(k_3)^2}{m}A  + (\lam L +  k_2)T e_{m-1}\right) + (\lam L+  k_2)\int_0^t \bfE \phi_m(x_{s^-})ds\\
&\leq& \left(\frac{C(k_3)^2}{m}A  + (\lam L +  k_2)T e_{m-1}\right)e^{(\lam L +  k_2)t}.
\eeao
Note that the same upper bound may be used in the case $g(x)>0$ since then $\phi_{m}(x_{s^-} +g(x_{s^-})) \leq \phi_{m}(x_{s^-}).$ 
Thus, using once more a property of $\phi$ and the above bound we have
$$
\bfE x^-_t - e_{m-1}\leq  \left(\frac{C(k_3)^2}{m}A  + (\lam L +  k_2)T e_{m-1}\right)e^{(\lam L +  k_2)t}
$$
and taking the limit $\mto$ implies $\bfE x^-_t\leq 0$ or $\bfP(x_t\leq0)=0$ for all $t\geq0.$
\epf

\bpf[Proof of Lemma \ref{CIRCEVdj:lem:inv_moment_bound}]
We apply the It\^o formula for the jump-extended CIR/CEV model and then take expectations to get\footnotesize
\beao
\bfE(x_t)^p &=& \bfE(\xi_0)^p + p\bfE\left(\int_0^t(x_{s^-})^{p-1}(k_1 - k_2x_{s^-}-\lam g(x_{s^-}))ds\right)\\
&& + \frac{(k_3)^2p(p-1)}{2}\bfE\left(\int_{0}^{t}(x_s)^{p-2+2\alpha}b(x_{s-\tau})^2ds\right) + \lam\bfE\left(\int_0^t \left((x_{s^-} +g(x_{s^-}))^p - (x_{s^-})^p \right)ds\right),
\eeao\normalsize
where $p<0.$  Now Assumption A yields $(x_{s^-} +g(x_{s^-}))^p \leq (x_{s^-})^p$ and we can find a constant $C>0$ such that 
$$
k_1p x^{p-1} +\frac{(k_3)^2p(p-1)}{2}x^{p-2+2\alpha}\leq C,
$$
thus
\beao
\bfE(x_t)^p &\leq& \bfE(\xi_0)^p  + C + CC_\gamma^2\int_0^t\bfE(x_{s-\tau})^{2\gamma}ds + C\int_0^t\bfE(x_{s^-})^{p}ds\\
&\leq& \bfE(\xi_0)^p + C + CC_\gamma^2At + C\int_0^t\bfE(x_{s^-})^{p}ds\\
&\leq& (\bfE(\xi_0)^p + C + CC_\gamma^2At)e^{Ct},
\eeao
by an application of the Gronwall inequality and Lemma \ref{CIRCEVdj:lem:moment_bound}, where the constant $C$ depends on $\bfE (|b(\xi_{-\tau})^2|).$ Therefore,
\beao
\sup_{-\tau\leq t\leq T}\bfE(x_t)^p &\leq& \sup_{\tau\leq t\leq 0}\bfE(x_t)^p + \sup_{0\leq t\leq T}\bfE(x_t)^p\\
&\leq&\bfE\|\xi \|^p + (\bfE(\xi_0)^p  + C + CC_\gamma^2At)e^{Ct}:=\hat{A}.
\eeao
\epf

\bpf[Proof of Theorem \ref{CIRCEVdj:thm:mean_reversionJASD}]
The proof of the mean reversion property of (JASDM) method requires moment bounds which are shown shortly after in Lemma \ref{CIRCEVdj-lemma:Lograte:SDuniformMomentBound}.
We take expectations in (\ref{CIRCEVdj:eq:JASD}), with a new constant $C$, to get, 
\beao
\bfE y_{t_{n+1}} &=& \bfE y_{t_{{n+1}^-}} + \bfE g(y_{t_{{n+1}^-}})\D \wt{N}_n\\
&\leq& \bfE y_{t_{{n+1}^-}} + \sqrt{\bfE g^2(y_{t_{{n+1}^-}})}\sqrt{\bfE(\D \wt{N}_n)^2}\\
&\leq& \bfE (z_{t_{n+1}})^2 + L\sqrt{\bfE (y_{t_{{n+1}^-}})^2}\sqrt{\lam \D}\leq \bfE (z_{t_{n+1}})^2 +  C\sqrt{\D}.
\eeao
Moreover\footnotesize 
\beao
\bfE (z_{t_{n+1}})^2 &\leq&  \bfE\left( y_{t_n}\left(1- \frac{k_2\D}{1+k_2\theta\D}\right) + \frac{k_1\D}{1+k_2\theta\D} - \frac{(k_3)^2}{4(1+k_2\theta\D)^2 }\frac{b^2(y_{t_{n-\tau}})}{(1+b(y_{t_{n-\tau}})\D^{m})^2}(y_{t_n})^{2\alpha-1}\D\right)\\
&& +\frac{(k_3)^2}{4(1+k_2\theta\D)^2} \bfE\left(\frac{b^2(y_{t_{n-\tau}})}{(1+b(y_{t_{n-\tau}})\D^{m})^2}(y_{t_n})^{2\alpha -1}(\D W_n)^2\right)\\
&\leq&  \left(1- \frac{k_2\D}{1+k_2\theta\D}\right)\bfE  y_{t_n} + \frac{k_1\D}{1+k_2\theta\D},
\eeao\normalsize
where we used that 
$$\bfE(\psi_{t_n}(\D W_n)^j) = \bfE\left(\bfE(\psi_{t_n}(\D W_n)^j| \bbf_{t_n})\right)=\bfE(\psi_{t_n}\bfE((\D W_n)^j| \bbf_{t_n})),$$
for $j=1,2$ and $\bfE (\D W_n)=0, \bfE (\D W_n)^2=\D$, which implies\footnotesize
$$
\bfE y_{t_{n+1}}\leq\left(1- \frac{k_2\D}{1+k_2\theta\D}\right)^n\bfE  \xi_0 + \left(\frac{k_1\D}{1+k_2\theta\D} +C\sqrt{\D}\right) \left(1- \frac{k_2\D}{1+k_2\theta\D}\right)^{n-1} + \frac{k_1\D}{1+k_2\theta\D} +C\sqrt{\D}.
$$\normalsize
Now, we take the limit $\nto$ when $\D(1-\theta)<\frac{1}{k_2}$ to get 
$\underset{\nto}{\lim} \bfE y_{t_n} \leq \frac{k_1}{k_2\theta}.$
\epf

\bpf[Proof of Theorem \ref{CIRCEVdj:thm:strong_convSDDE}]

We divide the proof is several steps. First we give moment bounds for the proposed numerical scheme. Then we estimate the error bound of the semi-discrete method. Using an auxiliary process $(h_t)$ we show strong convergence of $(h_t)$ to $(x_t^0)$ which carries on to the strong convergence of $(y_t^0)$ to $(x_t^0).$

\ble[Moment bounds for numerical approximation]\label{CIRCEVdj-lemma:Lograte:SDuniformMomentBound}
It holds that
$$
\bfE\sup_{\tau\leq t\leq T}(y_t)^p\leq A_p\bfE(\| \xi\| + k_1T)^p,
$$
for any $p>0,$ where $A_p$ is a constant.
\ele

\bpf[Proof of Lemma \ref{CIRCEVdj-lemma:Lograte:SDuniformMomentBound}]
We first observe that $(y_t)$ is bounded in the following way
\beao
0\leq y_t&\leq& \xi_0 + \int_0^t k_1ds + k_3\int_{0}^{t} \frac{b(y_{\hat{s}-\tau})}{1+b(y_{\hat{s}-\tau})\D^{m}}(y_{\hat{s}})^{\alpha-\frac{1}{2}}\sqrt{y_s}d\wt{W}_s\\
&\leq& \xi_0 + k_1T + k_3\int_{0}^{t} \frac{b(y_{\hat{s}-\tau})}{1+b(y_{\hat{s}-\tau})\D^{m}}(y_{\hat{s}})^{\alpha-\frac{1}{2}}\sqrt{y_s}d\wt{W}_s:=u_t,
\eeao
a.s., where the lower bound comes from the construction of $(y_t)$ and the upper bound follows from a comparison theorem. We will bound $(u_t)$ and therefore $(y_t),$ since $0\leq y_t\leq u_t$ a.s. Set the stopping time $\tau_R:=\inf\{t\in[0,T]: u_t>R\},$ for $R>0$ with the convention $\inf \emptyset=\infty.$ Application of It\^o's formula on $(u_{t\wedge\tau_R})^p$ implies\footnotesize
\beao
(u_{t\wedge\tau_R})^p&=&(\xi_0 + k_1T)^p + \frac{p(p-1)}{2}(k_3)^2\int_{0}^{t\wedge\tau_R} (u_s)^{p-2}\frac{b^2(y_{\hat{s}-\tau})}{(1+b(y_{\hat{s}-\tau})\D^{m})^2}(y_{\hat{s}})^{2\alpha-1}y_sds\\
&& + pk_3\int_{0}^{t\wedge\tau_R}(u_s)^{p-1}\frac{b(y_{\hat{s}-\tau})}{1+b(y_{\hat{s}-\tau})\D^{m}}(y_{\hat{s}})^{\alpha-\frac{1}{2}}\sqrt{y_s}d\wt{W}_s\\
&\leq&(\xi_0 + k_1T)^p + \frac{p(p-1)}{2}(k_3)^2\int_{0}^{t\wedge\tau_R} (u_s)^{p-1}b^2(y_{\hat{s}-\tau})(y_{\hat{s}})^{2\alpha-1}ds +
M_t\\
&\leq&(\xi_0 + k_1T)^p + \frac{p(p-1)}{2}(k_3)^2\int_{0}^{t\wedge\tau_R} \left(\frac{p-1}{2p}(u_s)^{p} + \frac{2^{p-1}b^{2p}(y_{\hat{s}-\tau})}{p}(y_{\hat{s}})^{(2\alpha-1)p}\right)ds +  M_t\\
&\leq&(\xi_0 + k_1T)^p + \frac{(p-1)^2}{4}(k_3)^2\int_{0}^{t\wedge\tau_R} (u_s)^{p}ds + 2^{3p-3}(p-1)(k_3)^2\int_{0}^{t\wedge\tau_R}b^{2p}(\xi_{-\tau})(y_{\hat{s}})^{(2\alpha-1)p}ds\\
&&+ 2^{3p-3}(p-1)(k_3C_{\gamma}^p)^2\int_{0}^{t\wedge\tau_R}
\left( (2\alpha-1)|y_{\hat{s}}|^p + 2(1-\alpha)|y_{\hat{s}-\tau}|^{p\gamma/(1-\alpha)}\right)ds +  M_t\\
&\leq&(\|\xi\| + k_1T)^p + \left(\frac{(p-1)^2}{4}+2^{3p-3}(p-1)[C_{\gamma}^p + b^{2p}(\xi_{-\tau})]\right) (k_3)^2\int_{0}^{t\wedge\tau_R} (u_s)^{p}ds\\
&&+ 2^{p-2}(p-1)2(1-\alpha)(k_3C_{\gamma}^p)^2\int_{0}^{t\wedge\tau_R}|u_{\hat{s}-\tau}|^{p\gamma/(1-\alpha)}ds +  M_t
\eeao\normalsize
where in the second step we have used that $0\leq y_t\leq u_t$ and $$M_t:=pk_3\int_{0}^{t\wedge\tau_R}(u_s)^{p-1}\frac{b(y_{\hat{s}-\tau})}{1+b(y_{\hat{s}-\tau})\D^{m}}(y_{\hat{s}})^{\alpha-\frac{1}{2}}\sqrt{y_s}d\wt{W}_s,$$ in the third step the inequality $x^{p-1}y\leq \ep\frac{p-1}{p}x^p + \frac{1}{p\ep^{p-1}}y^p,$ valid for $x\wedge y\geq0$ and $p>1$ with $\ep=\frac{1}{2}$ and in the final step the elementary inequality used in the proof of Lemma \ref{CIRCEVdj:lem:moment_bound} since  $\frac{1}{2}<\alpha<1.$ Taking expectations in the above inequality and using that $M_t$ is a local martingale vanishing at $0,$ we get
$$
\bfE(u_{t\wedge\tau_R})^p\leq \bba_p + \bbb_p\int_{0}^{t} \bfE(u_{s\wedge\tau_R-\tau})^{\frac{\gamma}{1-\alpha}p}ds + \bbc_p\int_{0}^{t} \bfE(u_{s\wedge\tau_R})^{p}ds,
$$
where the constants $\bba_p,\bbb_p$ and $\bbc_p$ depend on $p, \bfE\|\xi\|, k_1, k_2, k_3, T,C_\gamma$ and $\alpha.$ 
Now proceeding as in (\ref{CIRCEVdj:eq:mb4}) we get
$$
\bfE\sup_{-\tau\leq t\leq T}(u_{t})^{p}\leq A_p,
$$
and as a consequence
$$
\bfE\sup_{-\tau\leq t\leq T}(y_{t})^{p}\leq A_p,
$$
for any $p>0.$
\epf

\ble[Error bound]\label{CIRCEVdj-lemma:Lograte:SDErrorBounds}
Let $n_s$ be an  integer such that $s\in[t_{n_s},t_{n_s+1}].$ Then
$$
\bfE|y_s-y_{\hat{s}}|^p\leq\hat{A}_p\D^{p/2}, \qquad \bfE |y_s-y_{\wt{s}}|^p<\wt{A}_p\D^{p/2},
$$
for any $p>0,$ where the positive constants $\hat{A}_p,\wt{A}_p$ do not depend on $\D.$
\ele

\bpf[Proof of Lemma \ref{CIRCEVdj-lemma:Lograte:SDErrorBounds}]
First we take a $p\geq2.$ We get that
\footnotesize
\beao 
&&|y_s-y_{\hat{s}}|^p=\Big|\int_{t_{n_s}}^s\left(k_1 - k_2(1-\theta)y_{\hat{u}}-k_2\theta y_{\wt{u}}\right)du+ \int_{t_{n_s}}^{t_{n_s+1}}k_2\theta y_{\hat{s}}du  -\int_s^{t_{n_s+1}}k_2\theta y_{s}du\\
&&+ \int_s^{t_{n_s}}\left(k_1 - k_2(1-\theta)y_{t_{n_s}} - \frac{(k_3)^2}{4(1+k_2\theta\D)}\frac{b^2(y_{t_{n_s}-\tau})}{(1+b(y_{t_{n_s}-\tau})\D^{m})^2}(y_{t_{n_s}})^{2\alpha-1}\right)du\\
&& + k_3\int_{t_{n_s}}^{s} \frac{b(y_{\hat{u}-\tau})}{1+b(y_{\hat{u}-\tau})\D^{m}}(y_{\hat{u}})^{\alpha-\frac{1}{2}}\sqrt{y_u}d\wt{W}_u\Big|^p\\
&\leq&5^{p-1}\Big(\big|\int_{t_{n_s}}^s\left(k_1 - k_2(1-\theta)y_{\hat{u}}-k_2\theta y_{\wt{u}}\right)du\big|^p + (k_2)^p\theta^p(y_{\hat{s}})^p(t_{n_s+1}-t_{n_s})^p + (k_2)^p\theta^p(y_{s})^p(t_{n_s+1}-s)^p\\
&&+ \left|\int_s^{t_{n_s}}\left(k_1 - k_2(1-\theta)y_{t_{n_s}} - \frac{(k_3)^2}{4(1+k_2\theta\D)}\frac{b^2(y_{t_{n_s}-\tau})}{(1+b(y_{t_{n_s}-\tau})\D^{m})^2}(y_{t_{n_s}})^{2\alpha-1}\right)du\right|^p\\
&& + (k_3)^p\big|\int_{t_{n_s}}^{s}\frac{b(y_{\hat{u}-\tau})}{1+b(y_{\hat{u}-\tau})\D^{m}}(y_{\hat{u}})^{\alpha-\frac{1}{2}}\sqrt{y_u}d\wt{W}_u\big|^p\Big)\\
&\leq&5^{p-1}\Big(|t_{n_s}-s|^{p-1}\int_{t_{n_s}}^s\left|k_1 - k_2(1-\theta)y_{\hat{u}}-k_2\theta y_{\wt{u}}\right|^pdu + (k_2)^p\theta^p\left((y_{\hat{s}})^p+(y_{s})^p\right)\D^p\\
&&+ \left|k_1 - k_2(1-\theta)y_{t_{n_s}} - \frac{(k_3)^2}{4(1+k_2\theta\D)}\frac{b^2(y_{t_{n_s}-\tau})}{(1+b(y_{t_{n_s}-\tau})\D^{m})^2}(y_{t_{n_s}})^{2\alpha-1}\right|^p\D^p\\
&& + (k_3)^p\big|\int_{t_{n_s}}^{s} \frac{b(y_{\hat{u}-\tau})}{1+b(y_{\hat{u}-\tau})\D^{m}} (y_{\hat{u}})^{\alpha-\frac{1}{2}}\sqrt{y_u}d\wt{W}_u\big|^p\Big),
\eeao\normalsize
where we have used the Cauchy-Schwarz inequality. Taking expectations in the above inequality and using Lemma \ref{CIRCEVdj-lemma:Lograte:SDuniformMomentBound} and Doob's martingale inequality on the diffusion term we conclude
\beqq\label{CIRCEVdj-eq:Lograte:SDErrorBoundsp>2}
\bfE|y_s-y_{\hat{s}}|^p\leq\hat{A}_p\D^{p/2},
\eeqq
where the positive quantity $\hat{A}_p$ except on $p,$ depends also on the parameters $k_1,k_2,k_3,\theta, \alpha, \gamma$ but not on $\D.$ Now, for  $0<p<2$ we get
$$
\bfE|y_s-y_{\hat{s}}|^p\leq\left(\bfE|y_s-y_{\hat{s}}|^2\right)^{p/2}\leq\hat{A}_p\D^{p/2},
$$
where we have used Jensen's inequality for the concave function $\phi(x)=x^{p/2}.$ Following the same lines, we can show that
$$ 
\bfE|y_s-y_{\wt{s}}|^p\leq\wt{A}_p\D^{p/2},
$$ 
for any  $0<p,$ where the positive quantity $\wt{A}_p$ except on $p,$ depends also on the parameters $k_1,k_2,k_3,\theta, \alpha, \gamma$ but not on $\D.$
\epf

For the rest of the proof we rewrite the compact form of (\ref{CIRCEVdj-eq:SDsdeComp}) in the following way
\beqq\label{CIRCEVdj-eq:Lograte:SDcomp}
y_t=\un{\xi_0 + \int_0^t f_{\theta}(y_{\hat{s}},y_{\wt{s}})ds + \int_{0}^{t}g(y_{\hat{s}},y_s, y_{\hat{s}-\tau})d\wt{W}_s}_{h_t} + \int_t^{t_{n+1}}f_1(y_{t_n},y_{t}, y_{t_{n}-\tau})ds,
\eeqq
where $f_{\theta}(\cdot,\cdot)$ is given by (\ref{CIRCEVdj-eq:drift discetisation}). The auxiliary process $(h_t),$ with $h_t=\xi(t), t\in[-\tau,0]$ is close to $(y_t)$ as shown in the next result.

\ble[Moment bounds involving the auxiliary process]\label{CIRCEVdj-lemma:Lograte:MomentBoundsauxiliary}
For any $s\in[0,T]$ it holds that
\beqq\label{CIRCEVdj-eq:Lograte:MomentBoundsauxiliary}
\bfE|h_s-y_{s}|^p\leq C_p\D^{p},  \quad \bfE|h_s|^p\leq C_h,
\eeqq
and for $s\in[t_{n},t_{n+1}]$ we have that
$$ 
\bfE|h_s-y_{\hat{s}}|^p\leq \hat{C}_p\D^{p/2}, \qquad \bfE|h_s-y_{\wt{s}}|^p\leq \wt{C}_p\D^{p/2},
$$ 
for any $p>0,$ where  the positive quantities $C_p,\hat{C}_p,\wt{C}_p,C_h$ do not depend on $\D.$
\ele

\bpf[Proof of Lemma \ref{CIRCEVdj-lemma:Lograte:MomentBoundsauxiliary}]
We have that
$$
|h_s-y_s|^p=\left|\int_s^{t_{n+1}}f_1(y_{t_n},y_{s},y_{t_{n}-\tau})du \right|^p\leq|t_{n+1}-s|^p|f_1(y_{t_n},y_{s},y_{t_{n}-\tau})|^p,
$$
for any $p>0,$ where we have used (\ref{CIRCEVdj-eq:Lograte:SDcomp}). Using Lemma \ref{CIRCEVdj-lemma:Lograte:SDuniformMomentBound} we get the left part of (\ref{CIRCEVdj-eq:Lograte:MomentBoundsauxiliary}). Now for $p>2$ and noting that
\beao
\bfE|h_s|^p&\leq&2^{p-1}\bfE|h_s-y_s|^p + 2^{p-1}\bfE|y_s|^p\\
&\leq&2^{p-1}C_p\D^{p} + 2^{p-1}A_p\bfE(\| \xi \| + k_1T)^p\leq C_h,
\eeao
we get the right part of (\ref{CIRCEVdj-eq:Lograte:MomentBoundsauxiliary}), where we have used Lemma \ref{CIRCEVdj-lemma:Lograte:SDuniformMomentBound}. The case $0<p<2$ follows by Jensen's inequality as in Lemma \ref{CIRCEVdj-lemma:Lograte:SDErrorBounds}.

Furthermore, for $s\in[t_{n},t_{n+1}]$ and $p>2$ we derive that
\beao
\bfE|h_s-y_{\hat{s}}|^p&\leq&2^{p-1}\bfE|h_s-y_s|^p + 2^{p-1}\bfE|y_s-y_{\hat{s}}|^p\\
&\leq& 2^{p-1}C_p\D^{p} + 2^{p-1}\hat{A}_p\D^{p/2}\leq\hat{C}_p\D^{p/2},
\eeao
where we have used (\ref{CIRCEVdj-eq:Lograte:SDErrorBoundsp>2}) and in the same manner
$$
\bfE|h_s-y_{\wt{s}}|^p\leq 2^{p-1}C_p\D^{p} + 2^{p-1}\wt{A}_p\D^{p/2}\leq\wt{C}_p\D^{p/2}.$$
The case $0<p<2$ follows by Jensen's inequality.
\epf

The following results examine the convergence of the auxiliary process $(h_t)$ to $(x_t)$ in $\bbl^2$.

\bpr\label{CIRCEVdj-propo:L2auxiliaryConv}
Let Assumptions A and B hold and take $m=1/4$ in (\ref{CIRCEVdj-eq:SDsdeComp}). Then we have
\beqq \label{CIRCEVdj-eq:L2auxiliaryConv}
\bfE\sup_{-\tau\leq t\leq T}|h_{t}-x_{t}|^2\leq C\D^{(\alpha-\frac{1}{2})\wedge\gamma},
\eeqq
where $C$ is independent of $\D.$ In case $\alpha=1/2$
\beqq \label{CIRCEVdj-eq:L2auxiliaryConvCIR}
\bfE\sup_{-\tau\leq t\leq T}|h_{t}-x_{t}|^2\leq C\D^{\gamma}.
\eeqq
\epr

\bpf[Proof of Proposition \ref{CIRCEVdj-propo:L2auxiliaryConv}]
We estimate the difference $|\bbE_{t}|^2:=|h_{t}-x_{t}|^2.$ 

It holds that
\beam\nonumber
&&f_{\theta}(y_{\hat{s}},y_{\wt{s}})-f_{\theta}(x_s,x_s)=(k_1 - k_2(1-\theta)y_{\hat{s}}-k_2\theta y_{\wt{s}}) - (k_1-k_2x_s)\\
\nonumber&=& - k_2(1-\theta)(y_{\hat{s}}-x_s)-k_2\theta (y_{\wt{s}}-x_s)\\
\label{CIRCEVdj-eq:auxiliary_error_drift}&=&  k_2(1-\theta)(h_s - y_{\hat{s}}) + k_2\theta (h_s-y_{\wt{s}}) -k_2(h_s-x_s),
\eeam
therefore\footnotesize
\beao
&&|\bbE_{t}|^2=\left|\int_{0}^{t}\left(f_{\theta}(y_{\hat{s}},y_{\wt{s}})-f_{\theta}(x_s,x_s)\right)ds + \int_{0}^{t}\left(g(y_{\hat{s}},y_s,y_{s-\tau})\textup{sgn}(z_s)-k_3b(x_{s-\tau})x_s^\alpha\right)dW_s\right|^2\\
&\leq& 2T\int_{0}^{t}\left(k_2(1-\theta)|h_s - y_{\hat{s}}| + k_2\theta |h_s-y_{\wt{s}}| + k_2|\bbE_s|\right)^2ds + 2|M_t|^2\\
&\leq& 6T(k_2)^2(1-\theta)^2\int_{0}^{t}|h_s - y_{\hat{s}}|^2ds + 6T(k_2)^2\theta^2\int_{0}^{t}|h_s-y_{\wt{s}}|^2ds + 6T(k_2)^2\int_{0}^{t}|\bbE_s|^2ds + 2|M_t|^2,
\eeao\normalsize
where in the second step we have used the Cauchy-Schwarz inequality and (\ref{CIRCEVdj-eq:auxiliary_error_drift}) and
$$
M_t:= \int_{0}^{t} (g(y_{\hat{u}},y_{\wt{u}}, y_{\hat{u}-\tau})\textup{sgn}(z_u)-k_3b(x_{u-\tau})x_u^\alpha) dW_u.
$$
Taking the supremum over all $t\in[0,t_1]$ and then expectations we have\footnotesize
\beam
\nonumber
&&\bfE\sup_{0\leq t\leq t_1}|\bbE_t|^2\leq 6T(k_2)^2(1-\theta)^2\int_{0}^{t_1}\bfE|h_s - y_{\hat{s}}|^2ds +
6T(k_2)^2\theta^2\int_{0}^{t_1}\bfE|h_s-y_{\wt{s}}|^2ds\\
\nonumber&& + 6T(k_2)^2\int_{0}^{t_1}\bfE\sup_{0\leq l\leq s}|\bbE_{l}|^2ds + 2\bfE\sup_{0\leq t\leq t_1} |M_t|^2\\
\nonumber&\leq& 6T^2(k_2)^2(1-\theta)^2\hat{A}_2\D + 6T^2(k_2)^2\theta^2\wt{A}_2\D + 6T(k_2)^2\int_{0}^{t_1}\bfE\sup_{0\leq l\leq s}|\bbE_{l}|^2ds\\
\label{CIRCEVdj-eq:L2auxiliaryIneq}&&+ 8\bfE |M_{t_1}|^2,
\eeam\normalsize
where in the second  step we have used Lemma \ref{CIRCEVdj-lemma:Lograte:SDErrorBounds} and Doob's martingale inequality with $p=2,$ since  $M_t$ is an $\bbR$-valued martingale that belongs to $\bbl^2.$ Moreover,\footnotesize
\beao
&&|g(y_{\hat{s}},y_s,y_{\hat{s}-\tau})-g(x_s,x_s,x_{s-\tau})|^2=\left|k_3\frac{b(y_{\hat{s}-\tau})}{1+b(y_{\hat{s}-\tau})\D^{m}}(y_{\hat{s}})^{\alpha-\frac{1}{2}}\sqrt{y_s}-k_3\frac{b(x_{s-\tau})}{1+b(x_{s-\tau})\D^{m}}(x_s)^\alpha\right|^2\\
&\leq&(k_3)^2\left[\frac{b(y_{\hat{s}-\tau})}{1+b(y_{\hat{s}-\tau})\D^{m}}\left(\sqrt{y_s}\left((y_{\hat{s}})^{\alpha-\frac{1}{2}}-(y_{s})^{\alpha-\frac{1}{2}}\right)   + ((y_s)^{\alpha}-(x_s)^\alpha)\right)\right.\\
&&\left. + (x_s)^\alpha \left(\frac{b(y_{\hat{s}-\tau})}{1+b(y_{\hat{s}-\tau})\D^{m}} - \frac{b(x_{s-\tau})}{1+b(x_{s-\tau})\D^{m}} \right)\right]^2\\
&\leq&6(k_3)^2(C_\gamma^2 |y_{\hat{s}-\tau}|^{2\gamma} + b^{2}(\xi_{-\tau}))\left(y_s\left((y_{\hat{s}})^{\alpha-\frac{1}{2}}-(y_{s})^{\alpha-\frac{1}{2}}\right)^2  +\frac{ ((y_s)^\alpha-(x_s)^\alpha)^2((y_s)^{1-\alpha}+(x_s)^{1-\alpha})^2}{((y_s)^{1-\alpha}+(x_s)^{1-\alpha})^2}\right)\\
&& + 3(k_3)^2(x_s)^{2\alpha}(b(y_{\hat{s}-\tau})-b(x_{s-\tau}))^2 \\
&\leq&6(k_3)^2 \left( (C_\gamma^2|y_{\hat{s}-\tau}|^{2\gamma} + b^{2}(\xi_{-\tau}))\left(y_s|y_{\hat{s}}-y_{s}|^{2\alpha-1} + \frac{4\alpha^2|y_s-x_s|^2}{((y_s)^{1-\alpha}+(x_s)^{1-\alpha})^2}\right) + \frac{C_\gamma^2}(x_s)^{2\alpha}2|y_{\hat{s}-\tau}-x_{s-\tau}|^{2\gamma}\right)\\
&\leq&6(k_3)^2(C_\gamma^2 |y_{\hat{s}-\tau}|^{2\gamma} + b^{2}(\xi_{-\tau}))\left(y_s|y_{\hat{s}}-y_{s}|^{2\alpha-1} + 8\alpha^2\frac{|h_s-y_s|^2 + |h_s-x_s|^2}{((y_s)^{1-\alpha}+(x_s)^{1-\alpha})^2}\right)\\
&&+3(k_3C_\gamma)^22^{(2\gamma-1)^+}(x_s)^{2\alpha}(|h_{s-\tau}-y_{\hat{s}-\tau}|^{2\gamma} + |h_{s-\tau}-x_{s-\tau}|^{2\gamma}),
\eeao\normalsize
or\footnotesize
\beam\nonumber
&&|g(y_{\hat{s}},y_s,y_{\hat{s}-\tau})-g(x_s,x_s,x_{s-\tau})|^2\leq6(k_3)^2(C_\gamma^2 |y_{\hat{s}-\tau}|^{2\gamma} + b^{2}(\xi_{-\tau}))y_s|y_{\hat{s}}-y_{s}|^{2\alpha-1}\\
\nonumber&& + 6(k_3)^2(C_\gamma^2 |y_{\hat{s}-\tau}|^{2\gamma} + b^{2}(\xi_{-\tau}))8\alpha^2\frac{|h_s-y_s|^2 + |\bbE_s|^2}{((y_s)^{1-\alpha}+(x_s)^{1-\alpha})^2}\\
&&\label{CIRCEVdj-eq:auxiliary_error_diffusion3}
+ 3(k_3C_\gamma)^22^{(2\gamma-1)^+}(x_s)^{2\alpha} |h_{s-\tau} - y_{\hat{s}-\tau}|^{2\gamma}+ 3(k_3C_\gamma)^22^{(2\gamma-1)^+}(x_s)^{2\alpha} (\bbE_{s-\tau})^{2\gamma},
\eeam\normalsize
where we have used the inequality $(n^{\alpha}-l^{\alpha})(n^{1-\alpha}-l^{1-\alpha})\leq2\alpha|n-l|,$ that is true for all $n,l\geq0$ and $1/2\leq\alpha\leq1,$
and the property of H\"{o}lder continuous functions $|x^\alpha-y^\alpha|\leq|x-y|^\alpha$ for $\alpha\leq1.$ Furthermore,

\beam\nonumber
|g(x_s,x_s,x_{s-\tau})-k_3b(x_{s-\tau})x_s^\alpha|^2&\!\!\!\!\leq&\!\!\!(k_3)^2(x_s)^{2\alpha}b^2(x_{s-\tau})\left(\frac{b(x_{s-\tau})\D^{m}}{1+b(x_{s-\tau})\D^{m}}\right)^2\\
&\!\!\!\!\!\leq&\!\!\!\nonumber(k_3)^2(x_s)^{2\alpha}b^4(x_{s-\tau})\D^{2m}\\
\label{CIRCEVdj-eq:auxiliary_error_diffusion4}
&\!\!\!\!\!\leq&\!\!\!\!\!\!8(k_3)^2(x_s)^{2\alpha}(C_\gamma^4 |x_{s-\tau}|^{4\gamma} + b^{4}(\xi_{-\tau}))\sqrt{\D}, 
\eeam
where we chose $m=1/4.$ 

At this point, we use a stochastic time change. We define the process
\beqq\label{CIRCEVdj-eq:time_change}
\zeta(t):=\int_0^t \frac{1536 (k_3)^2\alpha^2(C_\gamma^2 |y_{\hat{s}-\tau}|^{2\gamma} + b^{2}(\xi_{-\tau}))}{\left[(y_s)^{1-\alpha} + (x_s)^{1-\alpha}\right]^2}ds
\eeqq
and the stopping time
$$
\tau_l:=\inf\{s\in[0,T]: 6T (k_2)^2s + \zeta(s)\geq l\}.
$$
The process $\zeta(t)$ is well defined since $x_t>0$ a.s. and $y_t\geq0.$ We find that\footnotesize
\beao
&&\!\!\!\!\!\!\!\!\!\!\!\!\!\!\!\!\bfE |M_{t_1}|^2=\bfE\left|\int_{0}^{t_1} |g(y_{\hat{s}},y_s, y_{\hat{s}-\tau})\textup{sgn}(z_s)-k_3b(x_{s-\tau})x_s^\alpha|dW_s\right|^2= \bfE\left(\int_{0}^{t_1} |g(y_{\hat{s}},y_s, y_{\hat{s}-\tau})\textup{sgn}(z_s)-k_3b(x_{s-\tau})x_s^\alpha|^2ds\right)\\
&\leq&\!\!\!\!\!24(k_3)^2\bfE\left(\int_{0}^{t_1} \left((C_\gamma^2 |y_{\hat{s}-\tau}|^{2\gamma} + b^{2}(\xi_{-\tau}))y_s|y_{\hat{s}}-y_{s}|^{2\alpha-1} +  (C_\gamma^2 |y_{\hat{s}-\tau}|^{2\gamma} + b^{2}(\xi_{-\tau}))8\alpha^2\frac{|h_s-y_s|^2 + |\bbE_s|^2}{((y_s)^{1-\alpha}+(x_s)^{1-\alpha})^2} \right)ds\right)\\
&&+12(k_3C_\gamma)^22^{(2\gamma-1)^+}\bfE\left(\int_{0}^{t_1}(x_s)^{2\alpha}|h_{s-\tau}-y_{\hat{s}-\tau}|^{2\gamma}ds\right) + 12(k_3C_\gamma)^22^{(2\gamma-1)^+}\bfE\left(\int_{0}^{t_1}(x_s)^{2\alpha}|\bbE_{s-\tau}|^{2\gamma}ds\right)\\
&&+ 32(k_3)^2\sqrt{\D}\bfE\left(\int_{0}^{t_1}(x_s)^{2\alpha}(C_\gamma^4 |x_{s-\tau}|^{4\gamma} + b^{4}(\xi_{-\tau}))ds\right)  +  2\bfE\left(\int_{0}^{t_1} g^2(y_{\hat{s}},y_s, y_{\hat{s}-\tau})(\textup{sgn}(z_s)-1)^2ds\right)\\
&\leq& 24(k_3)^2\int_{0}^{t_1}\sqrt{\bfE((C_\gamma^2 |y_{\hat{s}-\tau}|^{2\gamma} + b^{2}(\xi_{-\tau}))^2y_s^2)}\sqrt{\bfE(|y_{\hat{s}}-y_{s}|^{4\alpha-2})}ds + \frac{1}{8}\int_{0}^{t_1}(|h_s-y_s|^2 + |\bbE_s|^2)(\zeta_s)^{\prime}ds\\
&&+12(k_3C_\gamma)^22^{(2\gamma-1)^+}\int_{0}^{t_1}\sqrt{\bfE(x_s)^{4\alpha}}\sqrt{\bfE|h_{s-\tau}-y_{\hat{s}-\tau}|^{4\gamma}}ds+ 12(k_3C_\gamma)^22^{(2\gamma-1)^+}\int_{0}^{t_1}\bfE(x_s)^{2\alpha}|\bbE_{s-\tau}|^{2\gamma}ds\\
&&+ 32(k_3)^2\sqrt{\D}\int_{0}^{t_1}\sqrt{\bfE(x_s)^{4\alpha}}\sqrt{\bfE(C_\gamma^4 |x_{s-\tau}|^{4\gamma} + b^{4}(\xi_{-\tau}))^2}ds +2 \int_{0}^{t_1} \bfE g^2(y_{\hat{s}},y_s, y_{\hat{s}-\tau})\bbi_{\{z_s\leq0\}}ds,
\eeao\normalsize
where we have used (\ref{CIRCEVdj-eq:auxiliary_error_diffusion3}), (\ref{CIRCEVdj-eq:auxiliary_error_diffusion4}) and (\ref{CIRCEVdj-eq:time_change}). Now, Lemmata \ref{CIRCEVdj-lemma:Lograte:SDuniformMomentBound}, \ref{CIRCEVdj-lemma:Lograte:SDErrorBounds} and \ref{CIRCEVdj-lemma:Lograte:MomentBoundsauxiliary} and the fact that $\bfE g^2(y_{\hat{s}},y_s, y_{\hat{s}-\tau})\bbi_{\{z_s\leq0\}}\leq C\sqrt{\D}$ (see \cite[Lemma 3.2]{halidias:2015}, \cite[Section 5]{halidias_stamatiou:2015}) imply\footnotesize
\beao
\bfE |M_{t_1}|^2&\leq& C\sqrt{\D^{2\alpha-1}} + 
\frac{1}{8}\int_{0}^{\tau}\bfE |h_s - y_s|^2(\zeta_s)^{\prime}ds +  \frac{1}{8}\int_{0}^{\tau}\bfE|\bbE_s|^2(\zeta_s)^{\prime}ds
+ C\D^{\gamma}\\
&& + C\int_{0}^{t_1}\sqrt{\bfE(\bbE_{s-\tau})^{4\gamma}}ds + C\sqrt{\D}\\
&\leq& C\D^{(\alpha-\frac{1}{2})\wedge\gamma} + C\int_{0}^{t_1}\sqrt{\bfE |h_s - y_s|^4}\sqrt{\bfE((\zeta_s)^{\prime})^2}ds +  \frac{1}{8}\int_{0}^{t_1}\bfE|\bbE_s|^2(\zeta_s)^{\prime}ds + C\int_{0}^{t_1}\sqrt{\bfE(\bbE_{s-\tau})^{4\gamma}}ds\\
&\leq& C\D^{(\alpha-\frac{1}{2})\wedge\gamma} + C\D^2\int_{0}^{t_1}\sqrt{\bfE (x_s)^{2(1-\alpha)}}ds +  \frac{1}{8}\int_{0}^{t_1}\bfE|\bbE_s|^2(\zeta_s)^{\prime}ds + C\int_{0}^{t_1}\sqrt{\bfE(\bbE_{s-\tau})^{4\gamma}}ds\\
&\leq& C\D^{(\alpha-\frac{1}{2})\wedge\gamma}  +  \frac{1}{8}\int_{0}^{t_1}\bfE|\bbE_s|^2(\zeta_s)^{\prime}ds + C\int_{0}^{t_1}\sqrt{\bfE(\bbE_{s-\tau})^{4\gamma}}ds,
\eeao\normalsize
where we have used Lemma \ref{CIRCEVdj:lem:inv_moment_bound} and the asymptotic relations, $\D^l=o(\D^{\alpha-\frac{1}{2}})$ for all $l\geq\frac{1}{2}$ as $\D\downarrow0.$

Turning back to (\ref{CIRCEVdj-eq:L2auxiliaryIneq}) we get

\beqq\label{CIRCEVdj-eq:Pol:unifBound until stopping_time}
\bfE\sup_{0\leq t\leq t_1}|\bbE_t|^2\leq C\D^{(\alpha-\frac{1}{2})\wedge\gamma}  + \int_{0}^{t_1}\bfE\sup_{0\leq l\leq s}|\bbE_{l}|^2(6T(k_2)^2s + \zeta_s)^{\prime}ds+ C\int_{0}^{t_1}\sqrt{\bfE(\bbE_{s-\tau})^{4\gamma}}ds,
\eeqq\\
Relation (\ref{CIRCEVdj-eq:Pol:unifBound until stopping_time}) for $t_1=\tau\wedge\tau_l$ implies\footnotesize
\beam\nonumber
\bfE \sup_{0\leq t\leq \tau\wedge\tau_l}(\bbE_t)^2 &\leq& C\D^{(\alpha-\frac{1}{2})\wedge\gamma}  + \int_0^{\tau\wedge\tau_l} \bfE\sup_{0\leq r\leq s}(\bbE_r)^2 (6T (k_2)^2s + \zeta_s)^{\prime}ds + C\int_{0}^{\tau\wedge\tau_l}\sqrt{\bfE(\bbE_{s-\tau})^{4\gamma}}ds\\
\nonumber&\leq& C\D^{(\alpha-\frac{1}{2})\wedge\gamma}  + \int_0^l \bfE\sup_{0\leq j\leq u}(\bbE_{\tau_j})^2du\\
\label{CIRCEVdj-eq:Pol:unifBound until tau} &\leq& C\D^{(\alpha-\frac{1}{2})\wedge\gamma}e^l,
\eeam\normalsize
where in the last step we have used Gronwall's inequality. Using again relation (\ref{CIRCEVdj-eq:Pol:unifBound until stopping_time}) for $t_1\in[0,\tau]$ and under the change of variables $u=6T(k_2)^2s + \zeta_s$ we get\footnotesize
\beao
&&\bfE \sup_{0\leq t\leq t_1}(\bbE_t)^2\leq C\D^{(\alpha-\frac{1}{2})\wedge\gamma} + \int_0^{6(k_2)^2\tau T + \zeta_\tau}\bfE\sup_{0\leq l\leq \tau_u}(\bbE_l)^2du \\
&\leq&C\D^{(\alpha-\frac{1}{2})\wedge\gamma} + \int_0^{\infty}\bfE\left(\sup_{0\leq l\leq \tau_u}(\bbi_{\{6(k_2)^2\tau T + \zeta_\tau\geq u\}}\bbE_l)^2\right)du \\
&\leq&\int_0^{6(k_2)^2T^2}\bfE\sup_{0\leq l\leq \tau_u}(\bbE_l)^2du + \int_{6(k_2)^2T^2}^{\infty}\bfP(6(k_2)^2T^2 + \zeta_\tau\geq u)\bfE\left(\sup_{0\leq l\leq \tau_u}(\bbE_l)^2 \big| \{6(k_2)^2T^2 +\zeta_\tau\geq u\}\right)du\\
&&+ C\D^{(\alpha-\frac{1}{2})\wedge\gamma}\\
&\leq&C\D^{(\alpha-\frac{1}{2})\wedge\gamma}(e^{6(k_2)^2T^2}+1) +C\D^{(\alpha-\frac{1}{2})\wedge\gamma}\int_0^{\infty}\bfP(\zeta_\tau\geq u)e^udu,
\eeao\normalsize
where in the last steps we have used (\ref{CIRCEVdj-eq:Pol:unifBound until tau}). We proceed by showing that $u\rightarrow \bfP(\zeta_\tau\geq u)e^u\in\bbl^1(\bbR_+).$ Markov's inequality implies 
$$ 
\bfP(\zeta_\tau\geq u)\leq e^{-\ep u}\bfE(e^{\ep \zeta_\tau}),
$$ 
for any $\ep>0.$ Using (\ref{CIRCEVdj-eq:time_change}) and the fact that $x_s>0$ and $y_s\geq0$ we can provide the following bound
\beao
\zeta_\tau&=&\int_0^\tau \frac{1536 (k_3)^2\alpha^2(C_\gamma^2 |y_{\hat{s}-\tau}|^{2\gamma} + b^{2}(\xi_{-\tau}))}{\left[(y_s)^{1-\alpha} + (x_s)^{1-\alpha}\right]^2}ds \leq 768 (k_3C_\gamma)^2\alpha^2\int_0^\tau (y_{\hat{s}-\tau})^{2\gamma}(x_s)^{2\alpha-2}ds\\
&\leq& 1536 (k_3)^2\alpha^2(C_\gamma^2\sup_{-\tau\leq t\leq0}(y_{t})^{2\gamma} + b^{2}(\xi_{-\tau})) \int_0^\tau (x_s)^{2\alpha-2}ds,
\eeao   
thus
\beqq \label{CIRCEVdj-eq:Pol:exponentialMoment}
\bfE(e^{\ep\zeta_\tau})\leq \bfE\left(e^{\ep 1536 (k_3)^2\alpha^2(C_\gamma^2\sup_{-\tau\leq t\leq0}(\xi_{t})^{2\gamma} + b^{2}(\xi_{-\tau}))\int_0^\tau (x_s)^{2\alpha-2}ds} \right).
\eeqq

It remains to bound the exponential moments of $(x_t).$ We work as in \cite[Section 4]{halidias_stamatiou:2015}. First we find the dynamics of the transformation $v=x^{2-2\alpha}$ by application of It\^o's formula \footnotesize
\beao
v_t&=&v_0 + \int_0^t \left( \un{(1-2\alpha)(1-\alpha)(k_3)^2b^2(x_{s-\tau})}_{K_0} + \un{2(1-\alpha)k_1}_{K_1}(v_s)^{\frac{1-2\alpha}{2-2\alpha}} - \un{2(1-\alpha)k_2}_{K_2}v_s \right)ds\\
&& + \int_{0}^{t}\un{2k_3(1-\alpha)b(x_{s-\tau})}_{K_3}\sqrt{v_s}dW_s,
\eeao \normalsize
for $t\in [0,t_1],$ where $v_0=(\xi_0)^{2-2\alpha}>0.$ Then, by a comparison theorem \cite[Prop. 5.2.18]{karatzas_shreve:1988} we obtain that $v_t\geq \zeta_t(q)>0$ a.s. or $(v_t)^{-1}\leq(\zeta_t(q))^{-1}$ a.s. where the process $(\zeta_t(q))$ reads
\beqq\label{CIRCEVdj-eq:Pol:transformed_auxiliary}
\zeta_t(q)=\zeta_0 + \int_0^t (q - \left( K_2 + \eta(q)\right)\zeta_s)ds + \int_{0}^{t}K_3\sqrt{\zeta_s}dW_s,
\eeqq
for $t\in [0,t_1]$ with $\zeta_0(q)=v_0$  and 
$$
\eta(q) = \frac{(2\alpha-1)(q -  K_0)^{\frac{1}{2\alpha-1}}}{(k_1)^{\frac{2-2\alpha}{2\alpha-1}}}.
$$

Process (\ref{CIRCEVdj-eq:Pol:transformed_auxiliary}) is a square root diffusion process and when $\frac{2q}{(K_3)^2}-1\geq0$ or
\beqq\label{CIRCEVdj-eq:Pol:transformed_auxiliary_parameter}
q\geq 2(1-\alpha)^2(k_3)^2b^2(\xi_{s-\tau}),
\eeqq
remains positive if $\zeta_0(q)>0.$  Therefore,  $(x_t)^{2\alpha-2}\leq(\zeta_t(q))^{-1}$ a.s. so it suffices to bound exponential inverse moments of $(\zeta_t).$ By  \cite[Th. 3.1]{hurd_kuznetsov:2008} we have
\beqq\label{CIRCEVdj-eq:Pol:exponential_Inverse_Moment_CIR}
\bfE e^{\del\int_0^t(\zeta_s(q))^{-1}ds}\leq C_{HK}(\zeta_0)^{-\frac{1}{2}\left(\nu(q)-\sqrt{\nu(q)^2-8 \frac{\del}{(K_3)^2}}\right)},
\eeqq
for $0\leq \del\leq \left(\frac{2q}{(K_3)^2} - 1\right)^2\frac{(K_3)^2}{8}=:\nu(q)^2\frac{(K_3)^2}{8},$ where the positive constant $C_{HK}$ is explicitly given in \cite[(10)]{hurd_kuznetsov:2008} depends on the parameters $k_2,k_3,T,\alpha,$ but is independent of $\zeta_0.$ Thus the other condition that we require for parameter $q$ is
\beqq\label{CIRCEVdj-eq:Pol:exponential_Inverse_Moment_CIR_condition}
q\geq 2(1-\alpha)\sqrt{2\del}k_3b(\xi_{s-\tau}) + 2(1-\alpha)^2(k_3)^2b^2(\xi_{s-\tau}).
\eeqq
When (\ref{CIRCEVdj-eq:Pol:exponential_Inverse_Moment_CIR_condition}) is satisfied then (\ref{CIRCEVdj-eq:Pol:transformed_auxiliary_parameter}) is satisfied too, thus there is actually no restriction on the coefficient $\del$ in (\ref{CIRCEVdj-eq:Pol:exponential_Inverse_Moment_CIR}) since we can always choose appropriately a $q$ such that  (\ref{CIRCEVdj-eq:Pol:exponential_Inverse_Moment_CIR_condition})  holds satisfying $$q\geq2(1-\alpha)\sqrt{2\del}k_3\left(C_\gamma\|\xi\|^\gamma + b(0)\right)+ 2(1-\alpha)^2(k_3)^2\left(C_\gamma\|\xi\|^\gamma + b(0)\right)^2.$$ Relation (\ref{CIRCEVdj-eq:Pol:exponentialMoment}) becomes\footnotesize
\beqq\label{CIRCEVdj-eq:Pol:exponential_Inverse_Moment_CIR_bound}
\bfE(e^{\ep\gamma_\tau})\leq \bfE\left(e^{\ep 1536 (k_3)^2\alpha^2(C_\gamma^2\|\xi\|^{2\gamma} + b^{2}(\xi_{-\tau}))\int_0^{\tau} (v_s)^{-1}ds} \right)\leq \bfE\left(e^{\ep 1536 (k_3)^2\alpha^2(C_\gamma^2\|\xi\|^{2\gamma} + b^{2}(\xi_{-\tau}))\int_0^{\tau} (\zeta_s(q))^{-1}ds} \right).
\eeqq\normalsize
We therefore require that
\beqq\label{CEV-eq:Pol:paramaeter_condition}
1536 (k_3)^2\alpha^2(C_\gamma^2\|\xi\|^{2\gamma} + b^{2}(\xi_{-\tau}))\ep\leq \left(\nu(q)\right)^2\frac{(K_3)^2}{8}
\eeqq
and can always find a $\ep>1,$ such the above relation holds by choosing appropriately $q$ as discussed before. Relation (\ref{CIRCEVdj-eq:Pol:exponentialMoment}) becomes
$\bfE(e^{\ep\zeta_\tau})\leq A$
and therefore
$$
\bfP(\zeta_\tau\geq u)\leq Ae^{-\ep u},
$$
implying
\beao
\bfE \sup_{0\leq t\leq t_1}(\bbE_t)^2&\leq& C\D^{(\alpha-\frac{1}{2})\wedge\gamma}(e^{6(k_2)^2T^2}+1) +CA\D^{(\alpha-\frac{1}{2})\wedge\gamma}\int_0^{\infty}e^{(1-\ep)u}du\\
&\leq&C\D^{(\alpha-\frac{1}{2})\wedge\gamma},
\eeao
by choosing  $\ep>1.$ We apply again successively  (\ref{CIRCEVdj-eq:Pol:unifBound until stopping_time}) for $t_1=k\tau\wedge\tau_l$  for $k=2,\ldots,N_0$ 
to finally get 
\beqq\label{CIRCEVdj-eq:auxiliary_L2}
\bfE \sup_{0\leq t\leq T}(\bbE_t)^2\leq C\D^{(\alpha-\frac{1}{2})\wedge\gamma}.
\eeqq

For the CIR delay model, case $\alpha=1/2,$ we just sketch the main differences since we may follow the same steps. We now use the estimates, (see (\ref{CIRCEVdj-eq:auxiliary_error_diffusion3}) and (\ref{CIRCEVdj-eq:auxiliary_error_diffusion4})) 
\beam\nonumber
&&|g(y_s,y_{\hat{s}-\tau})-g(x_s,x_{s-\tau})|^2\leq 8(k_3)^2 (C_\gamma^2 |y_{\hat{s}-\tau}|^{2\gamma} +  b^{2}(\xi_{-\tau}))\frac{|h_s-y_s|^2 + |\bbE_s|^2}{(\sqrt{y_s} + \sqrt{x_s})^2}\\
&&\label{CIRCEVdj-eq:auxiliary_error_diffusion3CIR}
+ 2(k_3C_\gamma)^22^{(2\gamma-1)^+} x_s |h_{s-\tau} - y_{\hat{s}-\tau}|^{2\gamma}+ 2(k_3C_\gamma)^22^{(2\gamma-1)^+} x_s (\bbE_{s-\tau})^{2\gamma}
\eeam
and 
\beqq\label{CIRCEVdj-eq:auxiliary_error_diffusion4CIR}
|g(x_s,x_{s-\tau})-k_3b(x_{s-\tau})\sqrt{x_s}|^2 \leq 8(k_3)^2x_s( C_\gamma^2|x_{s-\tau}|^{4\gamma} +  b^{2}(\xi_{-\tau}))
\sqrt{\D},
\eeqq
where we chose again $m=1/4$ and, (see (\ref{CIRCEVdj-eq:diffusion discetisation}))
$$
g(y,z)=k_3 \frac{b(z)}{1+b(z)\D^{m}}\sqrt{y}.
$$
The process $\zeta$ reads
\beqq\label{CIRCEVdj-eq:time_changeCIR}
\zeta(t):=\int_0^t \frac{256(k_3)^2 (C_\gamma^2 |y_{\hat{s}-\tau}|^{2\gamma} + b^{2}(\xi_{-\tau}))}{(\sqrt{y_s} + \sqrt{x_s})^2}ds
\eeqq
and the stopping time now is
$$
\tau_l:=\inf\{s\in[0,T]: 4T (k_2)^2s + \zeta(s)\geq l\}.
$$
We get the following bound 
$$
\bfE \sup_{0\leq t\leq t_1}(\bbE_t)^2 \leq C\D^{\gamma\wedge(1/2)}(e^{4(k_2)^2T^2}+1) +C\D^{\gamma\wedge(1/2)}\int_0^{\infty}\bfP(\zeta_\tau\geq u)e^udu
$$
and it remains to bound the following exponential inverse moment, (see \ref{CIRCEVdj-eq:Pol:exponentialMoment})
\beqq \label{CIRCEVdj-eq:Pol:exponentialMomentCIR}
\bfE(e^{\ep\zeta_\tau})\leq \bfE\left(e^{\ep 8(k_3)^2 (C_\gamma^2 \|\xi\|^{2\gamma} +  b^{2}(\xi_{-\tau}))\int_0^\tau (x_s)^{-1}ds} \right).
\eeqq
The exponential inverse moment of the delay CIR model is finite, $\bfE e^{\del\int_0^t(x_s)^{-1}ds}<A,$  thus we can  find an $\ep>1$ such that (\ref{CIRCEVdj-eq:Pol:exponentialMomentCIR}) holds. We conclude as before,
\beqq\label{CIRCEVdj-eq:auxiliary_L2CIR}
\bfE \sup_{0\leq t\leq T}(\bbE_t)^2\leq C\D^{\gamma\wedge(1/2)}.
\eeqq
\epf

In order to finish the proof of Theorem \ref{CIRCEVdj:thm:strong_convSDDE} we just use the triangle inequality, Lemma \ref{CIRCEVdj-lemma:Lograte:MomentBoundsauxiliary}, Proposition \ref{CIRCEVdj-propo:L2auxiliaryConv} and (\ref{CIRCEVdj-eq:auxiliary_L2}) to get
\beao
\bfE\sup_{0\leq t\leq T}|y^0_{t}-x^0_{t}|^2 &\leq& 2\bfE\sup_{0\leq t\leq T}|h_{t}-y_{t}|^2 + 2\bfE\sup_{0\leq t\leq T}|\bbE_{t}|^2\\
&\leq& 2C\D^2 + 2C\D^{(\alpha-\frac{1}{2})\wedge\gamma}\leq C\D^{(\alpha-\frac{1}{2})\wedge\gamma},
\eeao
for the delay CEV jump model and (\ref{CIRCEVdj-eq:auxiliary_L2CIR}) for the CIR jump model
$$
\bfE\sup_{0\leq t\leq T}|y^0_{t}-x^0_{t}|^2 \leq 2C\D^2 + 2C\D^{\gamma\wedge(1/2)}\leq C\D^{\gamma\wedge(1/2)},
$$
where $C$  is independent of $\D.$
\epf

\bpf[Proof of Theorem \ref{CIRCEVdj:thm:strong_convSDDEjump}]

First we have to ensure that 
Denote $\bbE_{k^-}:=y_{t_{k^-}}-x_{t_{k^-}}.$ It holds that
\beao
|\bbE_{k}| &=&\left|\bbE_{k^-} +  \left(g(y_{t_{k^-}})-g(x_{t_{k^-}})\right)\D\wt{N}_{k-1}\right|\\
&\leq&|\bbE_{k^-}| +  |g(y_{t_{k^-}})-g(x_{t_{k^-}}| |\D\wt{N}_{k-1}|\\
&\leq&|\bbE_{k^-}| +  L|\bbE_{k^-}|(1+\lam\D)\\
&\leq&   \left(1+L(1+\lam\D)\right)|\bbE_{k^-}|
\eeao
where we used Assumption A. Therefore 
\beao
\bfE\sup_{1\leq k\leq N_T}|\bbE_k|^2&\leq&\left(1+L(1+\lam\D)\right)^2|\bbE_{k^-}|^2\\
&\leq&\left(1+L(1+\lam\D)\right)^2C\D^{(\alpha-\frac{1}{2})\wedge\gamma}\leq C\D^{(\alpha-\frac{1}{2})\wedge\gamma},
\eeao
by an application of Theorem \ref{CIRCEVdj:thm:strong_convSDDE}, where $C$ is a positive constant independent of $\D.$ We conclude that the jump adapted semi-discrete method (\ref{CIRCEVdj:eq:JASD}) converges in the mean square sense to the true solution of the jump-extended CEV model (\ref{CIRCEVdj:eq:sddej}) with order of convergence at least $((\alpha-\frac{1}{2})\wedge\gamma)/2.$ Moreover, the jump adapted semi-discrete method (\ref{CIRCEVdj:eq:JASD}) converges in the mean square sense to the true solution of the jump-extended CIR model (\ref{CIRCEVdj:eq:sddej}) with order of convergence  at least $((1/2)\wedge\gamma)/2.$
\epf

\subsection*{Conclusion}

In this work we have considered models with applications in finance, described by SDDEs with jumps. We have somehow unified existing models in the general model (\ref{CIRCEVdj:eq:sddej}) for which we prove uniqueness, positivity, moment boundness and the mean reversion property of the solution process. We then proposed the JASDM method to numerically approximate the solution process of the delay CIR/CEV model with jumps. The proposed scheme is nonnegative, strongly convergent in the mean square sense to the exact solution of (\ref{CIRCEVdj:eq:sddej}), with bounded moments and an analogue of the mean reversion property. We intend to make more numerical experiments to illustrate the impact of the diffusion exponent $\alpha$ and the delay coefficient $\gamma$ on the convergence of the (JASDM).
\\
\textbf{Acknowledgments.} The author would like to thank the anonymous referees for their useful suggestions. They
have improved the quality of the paper considerably. 

\bibliographystyle{unsrt}\baselineskip12pt 
\bibliography{CIR_CEVdj}


  
\end{document}